\documentclass{article}
\usepackage{amsmath} 
\usepackage[standard,thmmarks,amsmath]{ntheorem}
\usepackage{graphics} 
\usepackage{algorithm}
\usepackage{algpseudocode}
\usepackage{subfigure}
\usepackage{booktabs}
\usepackage{graphicx}
\usepackage{url} 

\DeclareMathOperator{\tr}{tr} \DeclareMathOperator{\diag}{diag}
\DeclareMathOperator*{\Argmin}{Arg\,min} \DeclareMathOperator*{\argmin}{arg\,min}
 
\DeclareMathOperator{\rnd}{rnd}

\def\le{\leqslant}
\def\leq{\leqslant}
\def\ge{\geqslant}

\newtheorem{note}{Note}

\begin{document}
\textwidth=30pc

\title{Lower bound for the cost of connecting tree with given vertex degree sequence}


\author{{\sc Mikhail Goubko}$^*$,\\[2pt]
V.A.~Trapeznikov Institute of Control Sciences of RAS,\\
117997, Prosfoyuznaya, 65, Moscow, Russia\\
$^*$Corresponding author: mgoubko@mail.ru\\[2pt]
{\sc Alexander Kuznetsov}\\[2pt]
Voronezh State University, Voronezh, Russia\\{avkuz@bk.ru}\\[6pt]}

\maketitle

\begin{abstract}
{The optimal connecting network problem generalizes many models of structure
optimization known from the literature, including communication and transport
network topology design, graph cut and graph clustering, structure identification
 from data, etc. For the case of connecting trees with the given
sequence of vertex degrees the cost of the optimal tree is shown to be bounded from
below by the solution of a semidefinite optimization program with bilinear matrix
constraints, which is reduced to the solution of a series of convex programs with
linear matrix inequality constraints. The proposed lower bound estimate is used to
construct several heuristic algorithms and to evaluate their quality on a variety of
generated and real-life data sets.} {Optimal communication network, generalized
Wiener index, origin-destination matrix, semidefinite programming, quadratic matrix
inequality.}
\\
2000 Math Subject Classification: 05C05, 05C07, 05C12, 05C35, 05C50, 68R10, 90C06,
90C22, 90C35, 90C59, 94C15
\end{abstract}

\section{Introduction}

The shortcut \emph{network} is used below for a simple connected undirected graph
with labeled vertices. So, networks with different labeling are considered distinct.

Let us consider the fixed \emph{set of terminals} $V =\{1,...,n\}$ indexed from $1$
to $n$ and denote a collection of networks over vertex set $V$ with $\Omega(V)$. Let
us assume we are given a symmetric non-negative \emph{flow matrix}
$A=(\mu_{ij})_{i,j=1}^n$ (where $\mu_{ij}$ is an informational or material flow
between the $i$-th and the $j$-th terminal, and \emph{set of admissible networks}
$\Omega\subseteq\Omega(V)$ (e.g., the set of all trees or of all bipartite graphs of
order $n$, etc.).

The \emph{optimal connecting network} (OCN) problem is that of finding an admissible
network $G^*\in \Omega$ with the minimum weighted average distance between vertex
pairs. In the other words, network $G^*\in \Omega$ is a solution of OCN problem if
and only if $C_A(G^*)\le C_A(G)$ for all $G\in \Omega$, where
\begin{equation}\label{eq_C_def}
C_A(G):=\sum_{\{i,j\} \subset V}\mu_{ij}d_G(i,j)=\frac{1}{2}\tr D(G)A.
\end{equation}

Here $d_G(i,j)$ is \emph{distance} between the $i$-th and the $j$-th vertices in
graph $G$, and $D(G)=(d_G(i,j))_{i,j=1}^n$ is the distance matrix of graph $G\in
\Omega$.

This framework, being simplistic at the first glance, however, has many classical
problems of combinatory optimization as special cases. Considering specific flow
matrices, sets of admissible networks, and specifying a concrete notion of graph
distance (the \emph{shortest-path distance}, the \emph{resistance distance}, or some
\emph{weighted distance}) one can obtain a \emph{quadratic assignment problem}
(QAP), a \emph{graph cut} or \emph{clustering problem}, or a sort of a problem of
structure identification from data (see examples in Section \ref{sec_literature}).

In this article we study a special case of OCN problem, which encapsulates the
essence of many difficulties that arise in OCN search. We consider the admissible
set, which contains all trees with the given sequence of vertex degrees, and the
(most popular) concept of the shortest-path graph distance.

If flow matrix $A$ has rank one, i.e., it can be represented as an outer product
$A=\mu\mu^\top$, where $\mu$ is some non-negative \emph{sequence of vertex weights},
cost function (\ref{eq_C_def}) reduces to the \emph{weighted Wiener index}
$WI_\mu(G)=\mu^\top D(G)\mu$ and OCN problem reduces to the recently solved problem
of the Wiener index optimization over the set of trees with given vertex weight and
degree sequences. In \cite{goubko2016wiener} the optimal tree is efficiently
constructed with a modification of the famous Huffman algorithm for the optimal
prefix code \cite{huffman1952method}.

Below we approximate the general flow matrix $A$ by a rank-one matrix obtaining a
lower-bound estimate for the optimal connecting tree cost. Calculation of the
estimate reduces to the non-convex semidefinite program. We solve it iteratively
through a series of constrained convex semidefinite programs effectively calculated
with standard optimization tools (we used CVX package with SeDuMi solver). It takes
reasonable time to calculate the estimate on a PC for trees with several hundreds of
vertices.

The quality of the lower bound is evaluated on a number of generated flow matrices
with dimension from 10 to 1000 and on the selected real-life origin-destination
matrices with dimension varying from 12 to 300. High quality of the lower bound is
verified in many practical cases, although in general the quality crucially depends
on how accurately matrix $A$ can be approximated by a rank-one matrix.

\section{Literature}\label{sec_literature}

\subsection{Quadratic assignment problems}\label{sec_QAP}

OCN problem is closely related to many structure optimization problems studied in
the literature. If all networks in the set $\Omega$ of admissible networks are
isomorphic and differ only in the vertex labeling, the solution of OCN reduces to
the assignment of terminals to network vertices, and we obtain a classical
Koopmans-Beckmann's quadratic assignment problem (QAP) \cite{koopmans1957assignment}
$$\min_\pi\sum_{i,j=1}^n\mu_{ij}d_{\pi(i)\pi(j)}\text{, where }\pi\text{ is a permutation of }1,...,n.$$
QAP is well-known as one of the most difficult problems of combinatory optimization
\cite{burkard1998quadratic,loiola2007survey,burkard2013quadratic}. It has many
unsolved instances of the dimension less than a hundred and does not have lower
bounds of guaranteed quality.

\subsection{Graph partitioning}\label{sec_graph_part}

If, in addition, the considered topology is a balanced tree of diameter $4$ with
$K+1$ internal vertices and only flows between tree leaves are allowed, the model is
equivalent to the optimal graph $K$-partitioning problem. If function $\pi(i)$
assigns a cluster number $1,..., K$ to $i$-th terminal $i=1,...,n$, then the cost
function reduces to
$$C_A(\cdot) = \sum_{k=1}^K\sum_{i:\pi(i)=k}\left[2\sum_{j:\pi(j)\neq k}\mu_{ij} +\sum_{j:\pi(j)=k}\mu_{ij}\right]=$$
$$=\sum_{k=1}^K\sum_{i\in s_k}\sum_{j \notin s_k}\mu_{ij}+\mathrm{const} = Cut(s_1, ..., s_K),$$
where $s_k:=\{i:\pi(i)=k\}$. In a similar fashion, \emph{balanced graph cut}
problems \cite{hein2011beyond} are obtained as a very special case of OCN.

Although the set of trees with the given sequence of vertex degrees includes the
admissible sets of graph partitioning and balanced cut problems (and even of the QAP
over the tree topology), the framework studied in this article is not completely
equivalent to the above problems.

The wider set does not necessary results in the more complex problem (e.g., the
complete graph is an obvious solution of the OCN over the set of all graphs of the
fixed order). On the other hand, the problem studied in the present article can be
seen as a variation of the balanced hierarchical clusterization problem, when not
only terminals have to be optimally grouped into clusters, but clusters should also
be rationally arranged into a hierarchy.

Business process decomposition and work breakdown structure (WBS) construction
problems are among possible applications. In many classic notations (IDEF, Aris,
BPMN, UML Activity Diagrams, Event Process Chains, and others) a business process in
an organization is represented as a directed graph where vertices are elementary
operations (activities) and arcs are labeled with material or information flows
between activities. In the same manner, vertices in a project schedule network are
project operations (works), while arcs represent precedence relations between them.

\begin{figure}[!h]
\centering\includegraphics[width=3.5in]{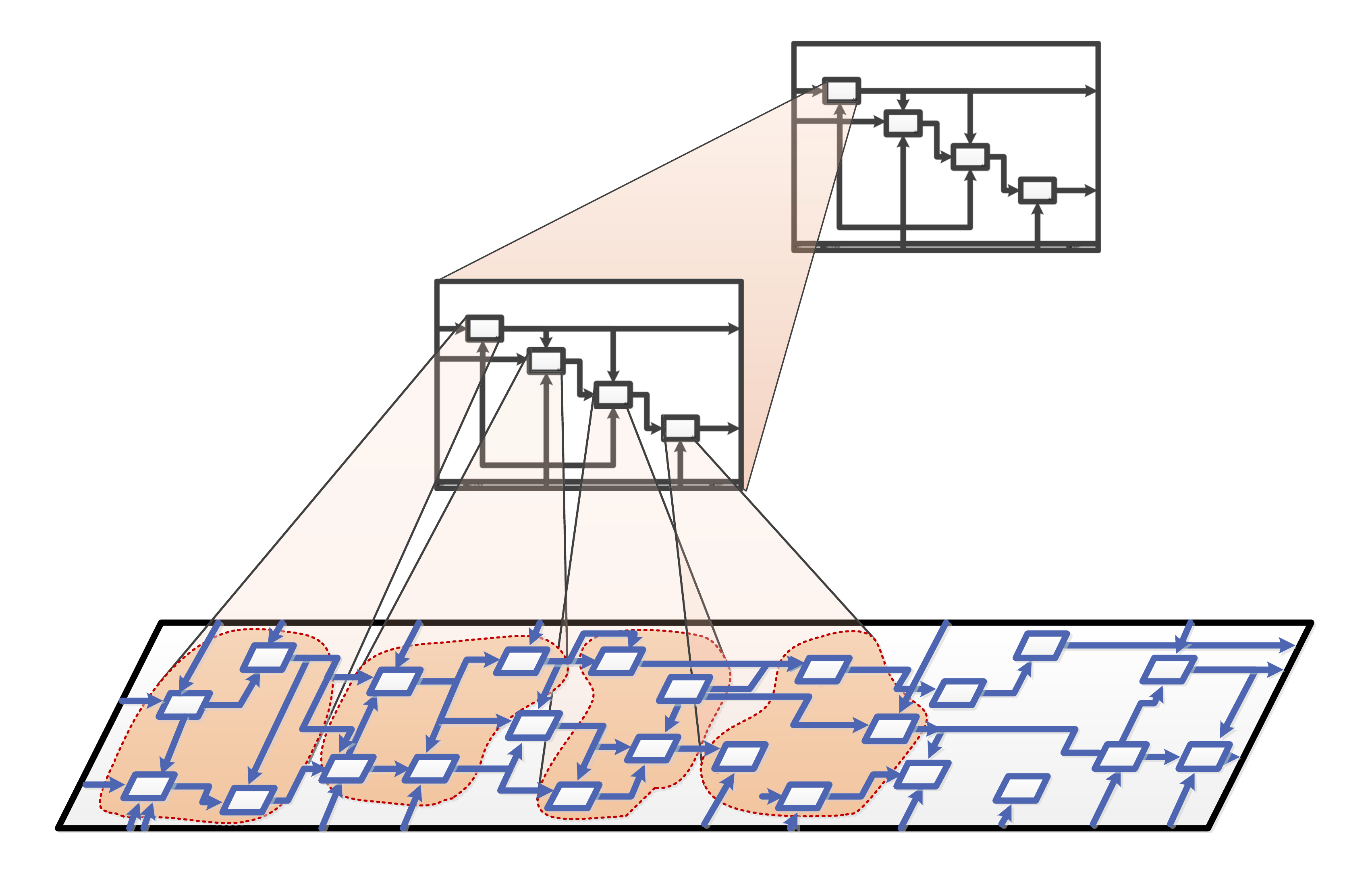}
\caption{Hierarchical decomposition of the business process} \label{fig_business}
\end{figure}

A complex business process (or a project schedule) may have many thousands
elementary activities. To simplify its representation and analysis, the activities
are arranged in a hierarchy of \emph{diagrams} so that only the limited number
(typically, from 5 to 7) of activities along with their internal and external flows
are combined in a single diagram (see Figure \ref{fig_business}) hiding the
complexity inside sub-diagrams.

During the business analysis \cite{caetano2012generation} most closely connected
activities (those having the maximum number of connecting flows or the maximum flow
volume between them) are located in a single diagram and are grouped together into a
corresponding combined activity. Then combined activities are grouped again at a
higher level of decomposition tree taking into account flows that connect them. It
is commonly recognized that such ``rational'' decomposition reveals the information
about internal structure of processes in an organization. In particular, business
process partitioning is used to identify services in SOA (service-oriented
architecture) \cite{ma2009evaluating}.

When a flow connects activities in different diagrams, it is depicted as an external
flow both in the source and in the destination diagrams. This flow is also copied as
external in all higher-level diagrams until the common parent diagram where it is
depicted as an internal flow (see Figure \ref{fig_business}).

Let the flow matrix $A$ be the adjacency matrix of the flow graph where the
direction of arcs is ignored, and the tree-shaped network $G$ coincide with the
process decomposition hierarchy. Then the total number of flows in all diagrams
(counting for flow copies in different diagrams) is given by expression
(\ref{eq_C_def}), and rational business process decomposition reduces to OCN over
the set of hierarchies (trees) with a limited maximum vertex degree (typically, it
varies from 6 to 8).

\subsection{Wiener index}\label{sec_wiener}

If $A$ is an all-ones matrix and $d_G(\cdot,\cdot)$ is the (edge) distance in graph
$G$, then $C_A(G)$ in (\ref{eq_C_def}) reduces to the \emph{sum of distances} in
graph $G$, also known as the \emph{Wiener index}, the one of the earliest and most
popular topological graph invariants widely used in mathematical chemistry and
network analysis as the measure of graph compactness. Compact connected graphs have
the small value of the Wiener index while more scattered graphs have the larger
index value. If $A = \mu\mu^\top$, where $\mu=(\mu_1,...,\mu_n)$ is a positive
sequence of vertex weights, $C_A(\cdot)$ becomes a variant of the Wiener index for
vertex-weighted graphs \cite{klavvzar1997wiener}.

Mathematical properties of the Wiener index and its extensions are studied for
decades by graph theorists (see the surveys in
\cite{dobrynin2001wiener,dobrynin2002wiener,aouchiche2006automated,gutman2010survey,gutman2013degree}).
Also, they also employed by many applications including mathematical chemistry
\cite{gutman2012distance}, analysis of social
\cite{freeman1978centrality,newman2001scientific} and communication
\cite{broder2000graph,imase1981design} networks.

Studies of extremal problems \cite{plesnik1984sum} is a valuable part the literature
on the Wiener index. In particular, Fischermann et al. \cite{fischermann2002wiener}
have shown a sort of balanced trees (aka Volkmann trees) to minimize the Wiener
index over the set of trees with the limited maximum vertex degree. The problem of
Wiener index minimization over the set of tress with the given degree sequence was
independently solved by \cite{zhang2008wiener,wang2008extremal} and the optimal tree
was characterized, being known as \emph{greedy tree} in \cite{wang2008extremal}, and
also as the \emph{breadth-first-search (BFS) tree} in \cite{zhang2008wiener}. Later
these results were extended to the Wiener index for vertex-weighted graphs. It has
been shown in \cite{goubko2016wiener,goubkomiloserdov2016wiener} that the,
so-called, \emph{generalized Huffman tree} minimizes the Wiener index over the set
of trees with given vertex weight and degree sequences. The present article is the
further extension of these results. Although no efficient exact solution is proposed
for the general flow matrix $A$, the cost of the generalized Huffman tree for the
conveniently chosen vertex weights gives the lower bound of the cost of the optimal
tree. Vertex weights corresponding to the best lower bound are calculated from a
non-convex optimization problem. They are also used in the heuristic algorithm to
efficiently construct a nearly optimal tree.

\subsection{Structure learning}\label{sec_struct_learn}

Another closely connected strand of the literature is learning the graph structure
from data. In the basic setting some signals (time series) are collected at the
vertices of an unknown graph and the problem is to elicit the edges (weighted, in
general), of the graph using correlation of signals in its vertices as a clue.
Typically, the lower the distance between signals is, the closer they should be
located in a graph.

Let $X=(x^{(1)}, ..., x^{(n)})^\top$ be an $m\times n$ matrix, where $x^{(i)}$ is an
$m$-dimensional row representing the signal located in vertex $i=1,...,n$ of an
unknown graph $G = \langle V,E\rangle$ with edge weights $w_{ij}$, $ij\in E$.

The search of the graph, in which $i$-th and $j$-th vertices are connected when the
distance $||x^{(i)}-x^{(j)}||$ between the corresponding signals is small, is often
(see \cite{kalofolias2016,egilmez2017graph} and the references therein) reduced to
the minimization of the function
\begin{equation}\label{eq_laplacian}\frac{1}{2}\sum_{i,j=1}^n w_{ij}||x^{(i)}-x^{(j)}||=\tr X^\top L(G)X=\tr
L(G)A,\end{equation} where $A:=XX^\top$ is the covariance matrix\footnote{It plays
the role of the flow matrix in these applications, so we use the same notation.} and
$L(G)$ is the Laplacian matrix of graph $G$:
$$L_{ij} = \left\{
\begin{array}{cc}
  -w_{ij} & i\neq j, \\
  \sum_k w_{ik} & i=j, \\
  0 & \text{otherwise}.\end{array}\right\}$$

The set of admissible graphs is additionally constrained to account for the \emph{a
priori} information about the target graph (e.g., maximum vertex degree,
connectedness, or edge density). Edge weights $w_{ij}$ are sought in
\cite{kalofolias2016} while in \cite{dong2016learning,egilmez2017graph} the authors
seek for the Laplacian matrix $L$ further relaxing the admissible set to the set of
all positive semidefinite matrices with zero row sums. Regularization terms are
added to (\ref{eq_laplacian}) in \cite{kalofolias2016,egilmez2017graph} to obey
local connectivity (every vertex must be connected to another vertex in a graph) and
obtain the desired graph density.

In the present article a similar problem is solved for the distance matrix on the
place of the graph Laplace matrix in (\ref{eq_laplacian}). Although both criteria
(\ref{eq_C_def}) and (\ref{eq_laplacian}) promote construction of the graph by
connecting vertices with highly correlated signals, their mathematical properties
are different. OCN is not directly reduced to the continuous (and even complex)
optimization problem as in \cite{dong2016learning,kalofolias2016,egilmez2017graph}.
Instead we construct a lower bound estimate using the OCN with the rank-one flow
matrix, for which an exact solution is known.


\section{Weighted Wiener index}

As noted in Section \ref{sec_QAP}, the general OCN is strongly NP-complete. At the
same time, efficient algorithms are known for special cases. For example, as soon as
the complete graph is admissible, it is an obvious solution of OCN problem.

The case of the flow matrix of rank one also appears computationally tractable. If
$A=\mu\mu^\top$, where $\mu_i\ge 0$ is a \emph{weight} of terminal $i=1,...,n$, then
$C_A(T)$ reduces to the vertex-weighted Wiener index $WI_\mu(T)=\mu^\top D(T)\mu$,
for which an optimal connecting tree for a given vertex degree sequence is
effectively built by the generalized Huffman algorithm \cite{goubko2016wiener}.

Below in this section we provide basic notation and definitions, and also introduce
the generalized Huffman algorithm, which is extensively used below.

Let $d_G(v)$ be the degree of vertex $v\in V$ in network $G\in \Omega(V)$.
\emph{Vertex degree sequence} of network $G$ is a vector $d_G = (d_G(i))_{i=1}^n$.
Vertex $v\in V$ is called \emph{pendent} if $d_G(v)=1$ and is called \emph{internal}
otherwise.
\begin{definition}
Connected network $T\in \Omega(V)$ is called a \emph{tree} if $\sum_{i=1}^n
d_T(i)=2(n-1)$. The collection of trees over vertex set $V$ is denoted with
$\mathcal{T}(V)$.
\end{definition}
\begin{definition}
Natural sequence $d=(d_1,...,d_n)$ is called \emph{generating} if $\sum_{i=1}^n
d_i=2(n-1)$. Let $\mathcal{T}(d):=\{T\in \mathcal{T}(V): d_T=d\}$ denote the
collection of trees with degree sequence $d$.
\end{definition}

Let $K_{W,M}$ be the \emph{complete bipartite network} over vertex subsets $W$ and
$M$, i.e., $K_{W,M}$ has vertex set $W \cup M$ and edge set $W \times M$.

For a fixed weight sequence $\mu\in \mathbb{R}^n_+$ and generating degree sequence
$d$ the \emph{generalized Huffman algorithm} \cite{goubkomiloserdov2016wiener}
builds a tree $H \in \mathcal{T}(d)$ as shown in Listing \ref{lst:huffman}.

\begin{algorithm}[!ht]
\floatname{algorithm}{Listing} \caption{Build a Huffman tree for weight sequence
$\mu$ and degree sequence $d$}\label{lst:huffman}
\begin{algorithmic}[1]

\Function{HuffmanTree}{$\mu$, $d$}

\State $V:=\{1,...,n\}$ \Comment{Vertex set: sequences $\mu$ and $d$ are assumed to
have $n$ components} \State $W := \{i\in V: d_i=1\}$ \Comment{Index set for vacant
pendent vertices} \State $M := V\backslash W$ \Comment{Index set for vacant internal
vertices} \State $H=\langle V,\emptyset\rangle$ \Comment{Start with empty network
over vertex set $V$}.

\For {$i=1$ to $q-1$} \State Choose any $m \in \Argmin\left\{d_u| u\in \Argmin_{v\in
M}\mu_v\right\}$\label{lst:huffman:choose_m} \Comment{$m$ has the least degree among
vertices} \State \Comment{of the least weight in $M$.} \For {$j=1$ to $d_m-1$}
\State Choose any $w_j \in \Argmin_{w\in W} \mu_w$\label{lst:huffman:choose_w}
\State $W := W \backslash \{w_j\}$ \State $\mu_m := \mu_m+\mu_{w_j}$ \EndFor
\Comment{Pick the vertices $w_1, ..., w_{d_m-1}$ that have $d_m-1$ least weights in
$W$.} \State $H := H \cup \{w_1m\} \cup ... \cup \{w_{d_m-1}m\}$. \Comment{Add edges
$\{w_1m\},...,\{w_{d_m-1}m\}$ to network $H$} \State $M := M \backslash \{m\}$, $W =
W \cup \{m\}$ \Comment{Move $m$ to index set $W$ of vacant pendent vertices}
 \EndFor
\State $H = H \cup K_{W,M}$ \Comment{Finish the Huffman tree by adding the star
$K_{W,M}$}\State \Return{$H$} \Comment{By construction, at this moment $|M|=1$, and
$d_m = |W|$, where $\{m\}= M$} \EndFunction
\end{algorithmic}
\end{algorithm}

\begin{note}
Like the ``classic'' Huffman algorithm, this algorithm requires $\mathcal{O}(n\ln
n)$ operations, and, so, is highly efficient.
\end{note}

\begin{note}
Some freedom of choice is allowed at lines \ref{lst:huffman:choose_m} and
\ref{lst:huffman:choose_w} of the algorithm, so, several distinct Huffman trees are
possible, all sharing the same value of $WI_\mu(\cdot)$. Let $\mathcal{H}(\mu,d)$ be
the collection of Huffman trees for weight sequence $\mu$ and degree sequence $d$.
\end{note}

\begin{definition}
Weights $\mu$ are \emph{monotone} in degrees $d$ if for all $i, j \in V$ from $1 <
d_i < d_j$ it follows that $\mu_i \le \mu_j$.\footnote{We omit here the technical
assumption $\mu_i>0 \Leftrightarrow d_i=1$ imposed in
\cite{goubkomiloserdov2016wiener} to simplify the proofs.}
\end{definition}

\begin{theorem}\label{theorem_huffman_vwwi}\cite{goubko2016wiener,goubkomiloserdov2016wiener} \sloppy If weights $\mu$ are
monotone in degrees $d$ and tree $T$ minimizes $WI_\mu(T)$ over $\mathcal{T}(d)$,
then $T\in \mathcal{H}(\mu,d)$ (i.e., $T$ is a Huffman tree).
\end{theorem}

\begin{note}
Huffman tree can be built for any weight sequence $\mu$ but Theorem
\ref{theorem_huffman_vwwi} may fail if weights are not monotone in degrees.
\end{note}

\begin{figure}
\subfigure[Huffman tree (vertex weights are shown in
circles)]{\includegraphics[width=0.50\textwidth]{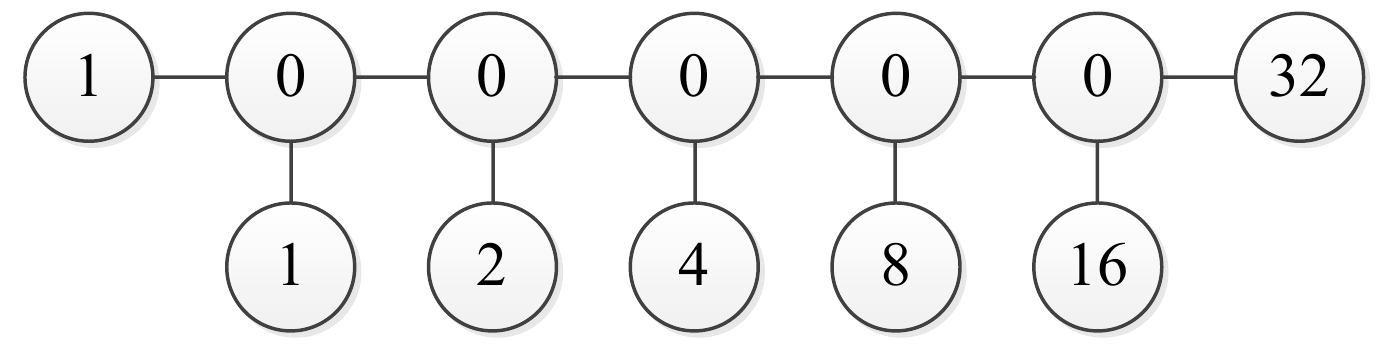}\label{fig_huffman}}
\subfigure[Greedy tree (\emph{aka}
BFS-tree)]{\includegraphics[width=0.40\textwidth]{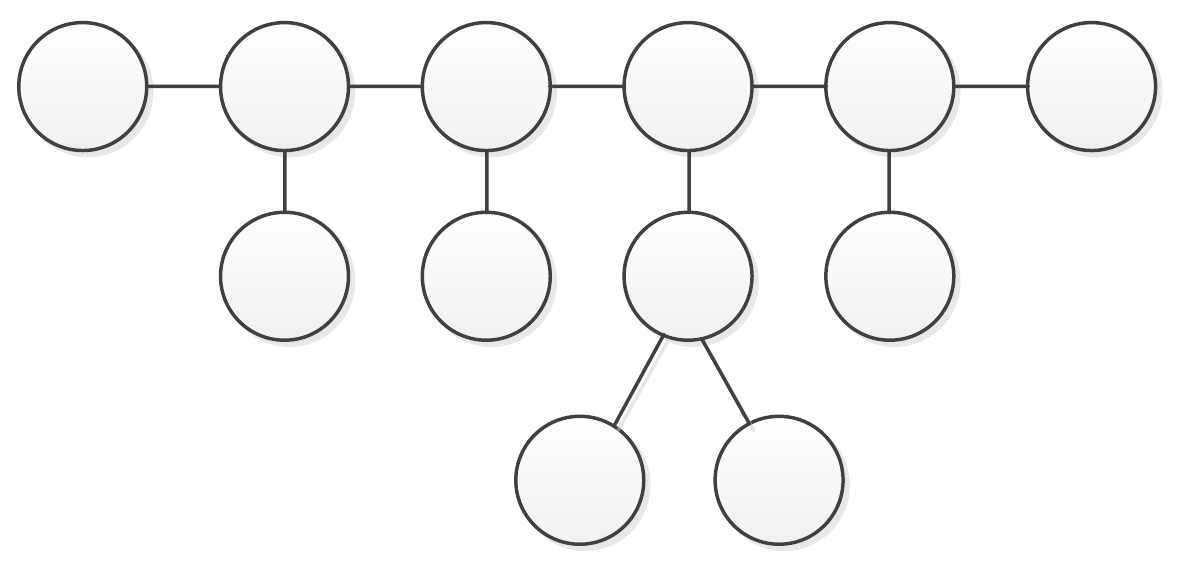}\label{fig_greedy}}
\caption{Huffman tree (a) not isomorphic to the greedy tree
(b)}\label{fig_huffman_greedy}
\end{figure}

\begin{note}
Only weights of internal vertices must be monotone in degrees in Theorem
\ref{theorem_huffman_vwwi}. Assume that, in addition, weights of pendent vertices
are required to not exceed those of internal vertices in $\mu$. Then, as shown in
\cite{cai2018ondistances}, all optimal trees for the degree sequence $d$ are
isomorphic to the \emph{greedy tree} (see Section \ref{sec_wiener} for details).
But, in general, Huffman trees may have diverse topology. For example, Huffman tree
for weight sequence $(1, 1, 2, 4, 8, 16, 32, 0, 0, 0, 0, 0)$ and degree sequence
$(1,1,1,1,1,1,1, 3, 3, 3, 3, 3)$ shown in Figure \ref{fig_huffman} is not isomorphic
to the corresponding greedy tree shown in Figure \ref{fig_greedy}.
\end{note}

\section{Lower bound of optimal connecting tree cost}\label{sec_lb}

In this article we study the following optimal connecting tree problem:
\begin{equation}\label{eq_min_problem}
\min_{T\in \mathcal{T}(d)} C_A(T) = \min_{T\in \mathcal{T}(d)} \tr D(T)A
\end{equation}
for given non-negative symmetric flow matrix $A$ and generating sequence
$d=(d_1,...,d_n)$ of vertex degrees, and in this section a closed-form expression is
derived for the lower bound estimate of the optimal tree cost. The main idea is to
approximate flow matrix $A$ by the sum of some non-negative rank-one matrix and a
diagonal matrix. The latter plays a role similar to that of the \emph{diagonal
perturbation} in \cite{rendl1995projection} and improves the quality of
approximation.

Let us denote $n\times n$ all-ones matrix with $J$ and define matrix $P(T) :=
\frac{n-1}{2} J-D(T)$. It is shown in \cite{bapat2005distance} that $P(T)$ is
positive semidefinite for any tree $T$ of order $n$.

\begin{theorem}\label{theorem_lb}
If real vector $\alpha \in \mathbb{R}^n$ and non-negative vector $\mu\in
\mathbb{R}^n_+$ are chosen such that weights $\mu$ are monotone in degrees $d$ and
matrix $\diag(\alpha)+\mu\mu^\top - A$ is positive semidefinite, then for any tree
$T\in \mathcal{T}(d)$
\begin{equation}\label{eq_def_LB}
C_A(T)\ge
LB(\alpha,\mu):=\frac{n-1}{2}\sum_{i,j=1}^n\mu_{ij}-\left(\frac{n-1}{2}\sum_{i=1}^n
\alpha_i+\mu^\top P(H(\mu))\mu\right),
\end{equation}
where $H(\mu)\in \mathcal{H}(\mu,d)$ is a Huffman tree for weight sequence $\mu$. In
other words, $LB(\alpha,\mu)$ is the \emph{lower bound estimate} for the problem
(\ref{eq_min_problem}).

\begin{proof}
Since matrices $P(T)$ and $\diag(\alpha)+\mu\mu^\top - A$ are positive semidefinite
and diagonal elements of $P(T)$ are equal to $\frac{n-1}{2}$,
\begin{multline}
C_A(T)=\tr D(T)A=\frac{n-1}{2}\sum_{i,j=1}^n\mu_{ij}-\tr P(T)A \ge\\
\ge\frac{n-1}{2}\sum_{i,j=1}^n\mu_{ij}-\frac{n-1}{2}\sum_{i=1}^n \alpha_i -\mu^\top
P(T)\mu =
\\
=\frac{n-1}{2}\sum_{i,j=1}^n(\mu_{ij}-\mu_i\mu_j)-\frac{n-1}{2}\sum_{i=1}^n
\alpha_i+2WI_\mu(T).
\end{multline}

From Theorem~\ref{theorem_huffman_vwwi} we know that $WI_\mu(T)\ge WI_\mu(H(\mu))$.
Hence,
\begin{multline}\label{eq_LB}
C_A(T)\ge \frac{n-1}{2}\sum_{i,j=1}^n\left(\mu_{ij}-\mu_i\mu_j\right)-\frac{n-1}{2}\sum_{i=1}^n \alpha_i+2WI_\mu(H(\mu))=\\
=\frac{n-1}{2}\sum_{i,j=1}^n\mu_{ij}-\left(\frac{n-1}{2}\sum_{i=1}^n
\alpha_i+\mu^\top P(H(\mu))\mu\right).
\end{multline}
This completes the proof.
\end{proof}
\end{theorem}

\section{Calculation of Lower bound}\label{sec_alg}

Inequality (\ref{eq_def_LB}) is valid for any combination of vectors $\alpha$ and
$\mu$ that satisfy conditions of Theorem \ref{theorem_lb}. Generally, we are
interested in the best (i.e., the largest) lower bound, which can be found by
maximizing $LB(\alpha,\mu)$ over all admissible combinations of $\alpha$ and $\mu$.
In this section we characterize the corresponding optimization problem, discuss its
algorithmic aspects and propose the optimization algorithm.

Taking into account Expression (\ref{eq_def_LB}), this problem is equivalent to the
minimization of the function
\begin{equation}\label{eq_min_criterion}
\mu^\top  P(H(\mu))\mu+\frac{n-1}{2}\sum_{i=1}^n \alpha_i.
\end{equation}
Since $\mu^\top P(G)\mu = \text{const} - 2WI_\mu(G)$ for fixed $\mu$ and any $G\in
\Omega(V)$, from Theorem \ref{theorem_huffman_vwwi} we know that
\begin{equation}\label{eq_Huffman_is_the_best}
\mu^\top  P(H(\mu))\mu = \max_{T\in \mathcal{T}(d)}\mu^\top  P(T)\mu,
\end{equation}
and, so, function  (\ref{eq_min_criterion}) is convex as an upper boundary of a
family of convex functions.

Finally, the best lower bound can be calculated from the minimization of a linear
function
\begin{equation}\label{eq_min_alpha_phi}
\min_{\alpha,\mu,\varphi} [\varphi+\frac{n-1}{2}\sum_{i=1}^n \alpha_i]
\end{equation}
under the \emph{bilinear matrix inequality} (BMI) constraint\footnote{Notation $
B\succeq 0$ means that matrix $B$ is positive semidefinite.}
\begin{equation}\label{eq_BMI}
\diag(\alpha)+\mu \mu^\top -A \succeq 0
\end{equation}
and convex constraints
\begin{gather}
\mu_i\ge 0, i=1,...,n,\label{eq_mu_nonnegative}\\
\varphi\ge \mu^\top  P(H)\mu \text{ for all }H\in \mathcal{H},\label{eq_sup_huffman}\\
\mu_j\ge\mu_i \text{ for all }
i,j:d_j-d_i=1,d_i>1,\label{eq_internal_weights_monotone}
\end{gather}
where $\mathcal{H}=\cup_{\mu}\mathcal{H}(\mu,d)$ is the collection of Huffman trees
for vertex degree sequence $d$ and all monotone weight sequences.

The number of trees in $\mathcal{H}$ is finite but large enough for the problem to
become intractable. At the same time, only the small number of inequalities in
(\ref{eq_sup_huffman}) are \emph{active} (i.e., make an equality at the optimal
point), which makes \emph{constraint generation} a promising idea.

Constraint generation is an approach to optimization problems with a large number of
constraints \cite{Ben-Ameur2006}. In our case it involves two steps that run in a
cycle. At the first step of iteration $t$ a relaxed problem (\ref{eq_min_alpha_phi})
containing only a subset $\mathcal{H}_t\subset \mathcal{H}$ of the constraints in
(\ref{eq_sup_huffman}) is solved. Then, at the second step, a special separation
procedure adds inequalities that are violated by the relaxed solution
$\alpha(t),\mu(t),\varphi(t)$ forming the set of constraints $\mathcal{H}_{t+1}$ for
the next iteration. The process is iterated until no violated inequality is found
(and, thus, $\mathcal{H}_t=\mathcal{H}_{t+1}$).

It is clear that if $\alpha(t),\mu(t),\varphi(t)$ is an optimal solution of the
relaxed problem (\ref{eq_min_alpha_phi}) for some constraint subset
$\mathcal{H}_t\subset \mathcal{H}$, and Huffman tree $H\left(\mu(t)\right)$ for
weight sequence $\mu(t)$ belongs to the set $\mathcal{H}_t$, then $\varphi(t)\ge
\mu(t)^\top H\mu(t)$ for any $H\in \mathcal{H}$, i.e., the relaxed solution is also
the optimal solution of problem (\ref{eq_min_alpha_phi}) with the complete
constraint set $\mathcal{H}$. On the contrary, if $H\left(\mu(t)\right) \notin
\mathcal{H}_t$, the relaxed solution cannot be the optimal solution for the complete
constraint set. Therefore, in our case the separation procedure just adds the tree
$H\left(\mu(t)\right)$ to the constraint set $\mathcal{H}_t$.

For the first iteration we take the constraint set $\mathcal{H}_1=\{H\}$ containing
only Huffman tree $H\in \mathcal{H}(\mathbf{1}, d)$ for all-ones weight sequence
$\mathbf{1}$ (\emph{aka} BFS-tree \cite{zhang2008wiener} \emph{aka} greedy tree
\cite{wang2008extremal}). Greedy tree is a good starting point because in Section
\ref{sec_numeric} it is shown that for large random flow matrices it is almost
always optimal. Numeric experiments also show that typically just a few constraint
generation iterations are enough to converge.

Unfortunately, even for the limited constraint set problem (\ref{eq_min_alpha_phi})
is not trivial, because BMI constraint (\ref{eq_BMI}) bounds a non-convex region due
to the bilinear term $\mu\mu^\top$ (mathematical properties of this region are
summarized in Appendix). At the same time, this BMI can be linearized with respect
to $\mu$ in the neighborhood of any point $\nu$ as follows. Inequality
(\ref{eq_BMI}) is equivalent to
$$\diag(\alpha)+\mu\nu^\top + \nu\mu^\top - \nu\nu^\top + (\mu-\nu)(\mu-\nu)^\top-A\succeq 0.$$
Suppressing the last term (which is an always non-negative and positive semidefinite
matrix) naturally gives the following linear matrix inequality (LMI) in $\alpha$ an
$\mu$:
\begin{equation}\label{eq_relax_LMI}
\diag(\alpha)+\mu\nu^\top + \nu\mu^\top - \nu\nu^\top-A\succeq0,
\end{equation}
which always bounds a convex region being a subset of the region bounded by BMI
(\ref{eq_BMI}).

Linearized problem (\ref{eq_min_alpha_phi}) with BMI (\ref{eq_BMI}) replaced with
LMI (\ref{eq_relax_LMI}) is a convex SDP (semidefinite program), which can be
conveniently coded using the disciplined programming notation of CVX package for
Matlab \cite{grant2008cvx} and efficiently solved by any available SDP solver like
SDPT4, SeDuMi, or Gurobi (we use SDPT4, the default solver for CVX shell).

To obtain the solution of the initial problem (\ref{eq_min_alpha_phi}) we combine
the majorization-minimization (MM) approach \cite{lange2016mm} with the alternating
directions (AD) method \cite{boyd2011distributed} solving in a cycle the linearized
problem and adjusting $\mu$ from the solution of the non-linearized problem under
fixed $\alpha$, the step, which is explained below.

Let us define symmetric matrix $A_\alpha:=A-\diag(\alpha)$, and denote its
eigenvalues $\lambda_i(A_\alpha)$ listed in the descending order, and the
corresponding eigenvectors $u^{(i)}(A_\alpha)$, $i=1,...,n$.

For fixed $\alpha$ BMI (\ref{eq_BMI}) is inconsistent whenever
$\lambda_2(A_\alpha)>0$ (see Lemma \ref{lemma_lambda2_empty}) and is satisfied for
any $\mu$ whenever $\lambda_1(A_\alpha)\le0$ (see Lemma
\ref{lemma_negative_semidefinite}). Otherwise (see Lemma \ref{lemma_analytical_eq}),
the region bounded by BMI is an interior of two convex sheets of a two-sheet
hyperboloid defined by the inequality
\begin{equation}\label{eq_hyper_quad}
\frac{\left(\mathbf{\mu}^\top u^{(1)}(A_\alpha)\right)^2}{\lambda_1(A_\alpha)}\ge
1-\sum_{i=2}^n\frac{\left(\mu^\top u^{(i)}(A_\alpha)\right)^2}{\lambda_i(A_\alpha)}.
\end{equation}
Alternatively the points satisfying (\ref{eq_hyper_quad}) are characterized by the
following pair (for ``$+$'' and for ``$-$'') of inequalities:
\begin{equation}\label{eq_hyper_conic}
\pm\frac{\mathbf{\mu}^\top u^{(1)}(A_\alpha)}{\sqrt{\lambda_1(A_\alpha)}}\ge
\sqrt{1+\sum_{i=2}^n\frac{\left(\mu^\top
u^{(i)}(A_\alpha)\right)^2}{|\lambda_i(A_\alpha)|}}.
\end{equation}
Absolute eigenvalues are used in (\ref{eq_hyper_conic}) to emphasize that
$\lambda_i(A_\alpha)\le0$ for all $i=2,...,n$. With notation
\begin{equation}\label{eq_hyper_vector}\mathbf{z}:=\left(1,\frac{\mu^\top u^{(2)}(A_\alpha)}{\sqrt{|\lambda_2(A_\alpha)|}},
..., \frac{\mu^\top
u^{(n)}(A_\alpha)}{\sqrt{|\lambda_n(A_\alpha)|}}\right)^\top\end{equation} conic
inequalities (\ref{eq_hyper_conic}) can be written in the canonic form
\begin{equation}\label{eq_conic_canonic}
\pm\frac{\mu^\top u^{(1)}(A_\alpha)}{\sqrt{\lambda_1(A_\alpha)}}\ge
\|\mathbf{z}\|_2.
\end{equation}

Therefore, for fixed $\alpha$, $\mu$-adjustment step reduces to the minimization of
$\varphi$ with respect to $\mu$ and $\varphi$ under constraints
(\ref{eq_mu_nonnegative}), (\ref{eq_sup_huffman}),
(\ref{eq_internal_weights_monotone}), and (\ref{eq_conic_canonic}) (for ``$+$'' and
for ``$-$''). This pair of conic programs is efficiently coded with CVX and solved
using almost any available convex programming tool (CPLEX, SDPT4, SeDuMi, Gurobi,
\emph{etc.}).\footnote{Due to nonnegativity and monotonicity constraints
(\ref{eq_mu_nonnegative}) and (\ref{eq_internal_weights_monotone}) one of these
programs is typically inconsistent, which does not make a problem.} Finally, the
adjusted $\mu$ is used as a new linearization point at the next iteration of the
algorithm.

\begin{note}
If $\lambda_i(A_\alpha)=0$ for some $i = 1,...,n$, the corresponding term in conic
inequalities (\ref{eq_hyper_quad})-(\ref{eq_conic_canonic}) is omitted, and the new
condition $\mu^\top u^{(i)}(A_\alpha)=0$ is added instead.
\end{note}

We need a feasible starting point to begin iterations. Lemma \ref{eq_X_is_empty}
says that for the feasible set to be not empty, $\alpha$ must be chosen such that
$\lambda_2(A_\alpha)\leq 0$. Therefore, let us choose
$\alpha(0)=\lambda_2(A)\mathbf{1}$, so that $\lambda_i(A_{\alpha(1)})=
\lambda_i(A)-\lambda_2(A)$. Eigenvectors of matrix $A_{\alpha(1)}$ coincide with
those of matrix $A$, so, according to Lemma \ref{lemma_perron_vector}, let us choose
feasible $\mu(1):=\sqrt{\lambda_1(A)-\lambda_2(A)}u^{(1)}(A)$, which can be used as
the first linearization point in (\ref{eq_relax_LMI}).

Function \textproc{MaximizeLB} that solves problem (\ref{eq_min_alpha_phi}) under
constraints (\ref{eq_BMI})-(\ref{eq_internal_weights_monotone}) is presented in
Listing \ref{lst:alg}. Combination of MM and AD steps highly improves convergence
compared to MM and AD applied separately.

\begin{algorithm}[!ht]
\floatname{algorithm}{Listing} \caption{Calculate the best parameters for the lower
bound $LB(\alpha,\mu)$}\label{lst:alg}
\begin{algorithmic}[1]

\Function{MaximizeLB}{} \State
$\mathcal{H}_1:=\left\{\text{\textproc{HuffmanTree}}(\mathbf{1},d)\right\}$
\Comment{Start from a single constraint in (\ref{eq_sup_huffman})} \State $t:=0$
\Repeat \State $t:=t+1$ \State $(\alpha(t), \mu(t), \varphi(t))
:=$\Call{SolveRelaxed}{$\mathcal{H}_t$} \State
$H:=\text{\textproc{HuffmanTree}}(\mu(t), d)$ \Comment{$H$ is a Huffman tree for
weight sequence $\mu(t)$} \State $\mathcal{H}_{t+1}=\mathcal{H}_t\cup\{H\}$
\Comment{Extend the set of constraints} \Until($\varphi(t)\ge\mu(t)^\top
P(H)\mu(t)$) \State \Return{$(\alpha(t),\mu(t))$} \EndFunction

\Function{SolveRelaxed}{$\mathcal{H}$} \State $t:=0$
\State $\mu(0):=\sqrt{\lambda_1(A)-\lambda_2(A)}u^{(1)}(A)$ \State $t:=0$ \Repeat
\State $t:=t+1$ \State $\alpha(t):=$\Call{SolveLinearized}{$\mathcal{H}, \mu(t-1)$}
\State $(\mu(t),\varphi(t)):=$\Call{AdjustMu}{$\mathcal{H}, \alpha(t)$}
\Until{$|\varphi(t-1)-\varphi(t)+\frac{n-1}{2}\sum_{i=1}^n
 [\alpha_i(t-1)-\alpha_i(t)]|<\delta$} \Comment{Improvement below tolerance} \State \Return{$(\alpha(t), \mu(t), \varphi(t))$}
\EndFunction

\Function{SolveLinearized}{$\mathcal{H}, \nu$}
 \State Find $(\alpha^*,\mu^*,\varphi^*)\in \Argmin_{\alpha,\mu,\varphi}[\varphi+\frac{n-1}{2}\sum_{i=1}^n
 \alpha_i]$ under constraints (\ref{eq_mu_nonnegative})-(\ref{eq_relax_LMI})\Comment{Convex SDP}
\State \Return{$\alpha^*$}
 \EndFunction

\Function{AdjustMu}{$\mathcal{H}, \alpha$}
 \State Find $(\mu^*,\varphi^*)\in \Argmin_{\mu,\varphi} \varphi$ under constraints (\ref{eq_mu_nonnegative})-(\ref{eq_internal_weights_monotone}),(\ref{eq_conic_canonic})\Comment{Pair of conic programs}
\State \Return{$(\mu^*,\varphi^*)$} \EndFunction
\end{algorithmic}
\end{algorithm}

\begin{note}Since the linearized solution is always feasible, $\mu^*$ and $\varphi^*$ from
\textproc{SolveLinearized} can be a starting point in \textproc{AdjustMu} for
algorithms that require an internal starting point.
\end{note}

\begin{note}
Two conic problems are solved in \textproc{AdjustMu}, one for ``$+$'' sign and the
other for ``$-$'' sign in (\ref{eq_conic_canonic}). However, for instance,
\texttt{cplexqcp} utility of CPLEX package solves both in a single run taking
inequality (\ref{eq_hyper_quad}) as an input.
\end{note}

\begin{note}The algorithm in \textproc{SolveRelaxed} converges,
since objective function (\ref{eq_min_alpha_phi}) is bounded, and every iteration
improves the solution. Numeric tests in Section \ref{sec_numeric} show fast
convergence in average (less than in a dozen iteration), however, in general, no
fast convergence can be guaranteed.
\end{note}

\begin{note}
To find the best values of the parameters of the lower bound the algorithm solves
the non-convex optimization problem. For such problems there is no universal
criterion of convergence to the global optimum. At the same time, global optimality
is not critical for lower bound evaluation since any admissible solution of problem
(\ref{eq_min_alpha_phi}) with constraints
(\ref{eq_BMI})-(\ref{eq_internal_weights_monotone}) gives a lower bound.
\end{note}

\section{Heuristics}\label{sec_heur}
One of applications of the lower bound estimate introduced in Section \ref{sec_lb}
is performance evaluation of heuristic algorithms that build nearly optimal trees
for the given degree sequence. Since any heuristic algorithm gives an upper bound to
the optimal tree cost, the gap between the upper and the lower bounds measures the
possible performance loss, justifies the price of algorithm improvement, and
motivates future research.

Heuristic algorithms may base on different ideas. In this section we describe two
algorithms that employ rank-one approximation of the flow matrix and the optimality
of Huffman trees.

Approximation of flow matrix $A$ with some matrix $\mu\mu^\top$ of rank one results
in assigning non-negative weights $\mu_i,i=1,...,n$, to the terminals. It is known
that the first principal component of a non-negative symmetric matrix $A$ is its
Perron vector $u^{(1)}(A)$. This means that $u^{(1)}(A) = \argmin_\mu
\|A-\mu\mu^\top\|_2$, so, the Perron vector is the best approximation (in $L_2$
norm) of matrix $A$ by a rank-one matrix. This justifies the choice of weight
sequence $u^{(1)}(A)$. By Perron-Frobenius theorem, the Perron vector is positive,
so,
 $u^{(1)}(A)$ is a valid weight sequence, to which the Huffman algorithm can be applied.
 Although $u^{(1)}(A)$ has not be monotone with respect to degree sequence $d$ and,
 so, Theorem \ref{theorem_huffman_vwwi} may not hold, the topology of the Huffman tree is still a
 good choice for a connecting tree with weight sequence $\mu$. Hence we introduce
 $$\text{\textproc{Heuristics1}}:=\text{\textproc{HuffmanTree}}(u^{(1)}(A), d).$$

Another low-rank approximation of the flow matrix goes from the lower bound
calculation (see the previous section). For
$(\cdot,\mu^{[2]})=\text{\textproc{MaximizeLB}}$ let us define
 $$\text{\textproc{Heuristics2}}:=\text{\textproc{HuffmanTree}}(\mu^{[2]}, d).$$
 The advantage of $\textproc{Heuristics2}$ is that weight sequence $\mu^{[2]}$ is always, by construction, monotone
 with respect to $d$ and, therefore, the Huffman tree is an optimal connecting tree
 for the approximated flow matrix. We postpone comparative performance analysis of
 both heuristics to the next section.

\section{Numeric simulations}\label{sec_numeric}

Several numeric tests on generated and real-world data were run to evaluate the
quality of the lower bound estimate proposed in Section \ref{sec_lb} compared to the
quality of two heuristic algorithms introduced in Section \ref{sec_heur}. The
performance is also estimated of the algorithm (see Section \ref{sec_alg}) for
calculation of the best parameter values of the lower bound.

\subsection{Random rank-one flow matrices}
First we check that the lower bound is tight when flow matrix $A$ has rank one. 100
degree sequences were generated for trees of order from 50 to 250 with degrees of
internal vertices uniformly distributed from 2 to 5. For every degree sequence $d$ a
monotone random vertex weight sequence $\nu$ was generated such that
$\nu_i=\rnd^\beta$, where $\rnd$ is a random number uniformly distributed on
$[0,1]$, and $\beta\ge0$ is a diversity factor (for $\beta=0$ all weights are equal
to unity, for $\beta=1$ we have the uniform distribution of weights, while for large
$\beta$ most weights, except some outliers, are close to zero).

The flow matrix $A$ was set to $\nu\nu^\top-\diag(\nu^2_1,..., \nu^2_n)$ (diagonal
entries of flow matrix are equal to zero). For all cases \textproc{MaximizeLB} was
called to find the best parameters of the lower bound. Two upper bounds and the
corresponding nearly optimal trees $H_1:=\text{\textproc{Heuristics1}}$ and
$H_2:=\text{\textproc{Heuristics2}}$ were obtained along with the
breadth-first-search tree $BFS := \text{\textproc{HuffmanTree}}(\mathbf{1},d)$. The
``best found tree'' was selected as $H^*=\Argmin_{H\in\{H_1,H_2,BFS\}}C_A(H)$.
Finally, we calculated the average cost $C_{avg}$ of 100 random trees from
$\mathcal{T}(d)$.

In all cases less than four iterations inside a single run of
\textproc{SolveRelaxed} function were enough to find the best parameters ($\mu =
\nu$, $\alpha=-\diag(\nu^2_1,..., \nu^2_n)$) that approximate perfectly this
simplistic flow matrix. So, it is no wonder that in all cases both heuristic
algorithms returned $T^*:=\text{\textproc{HuffmanTree}}(\nu)$, which was an exact
solution, since the lower bound $LB(\alpha,\mu)$ gave exactly $C_A(T^*)$. Thus, the
lower bound is tight in this setting.

$BFS$ tree is a ``perfectly balanced tree'' that can be calculated once for degree
sequence $d$ and used as a ``universal solution'' being more or less good for all
monotone vertex weight sequences. The relative gap $\Delta_{BFS}=
\frac{C_A(BFS)-C_A(T^*)}{C_A(T^*)}$ between the cost of $BFS$ and the cost of the
best found solution $H^*$ shows the price of knowing flow matrix $A$. The relative
gap $\Delta_{avg}= \frac{C_{avg}-C_A(T^*)}{C_A(T^*)}$ shows the price of solving OCN
problem in comparison with picking a random tree as a solution.

\begin{figure}\centering
\subfigure[Relative BFS gap
$\Delta_{BFS}$]{\includegraphics[width=0.80\textwidth]{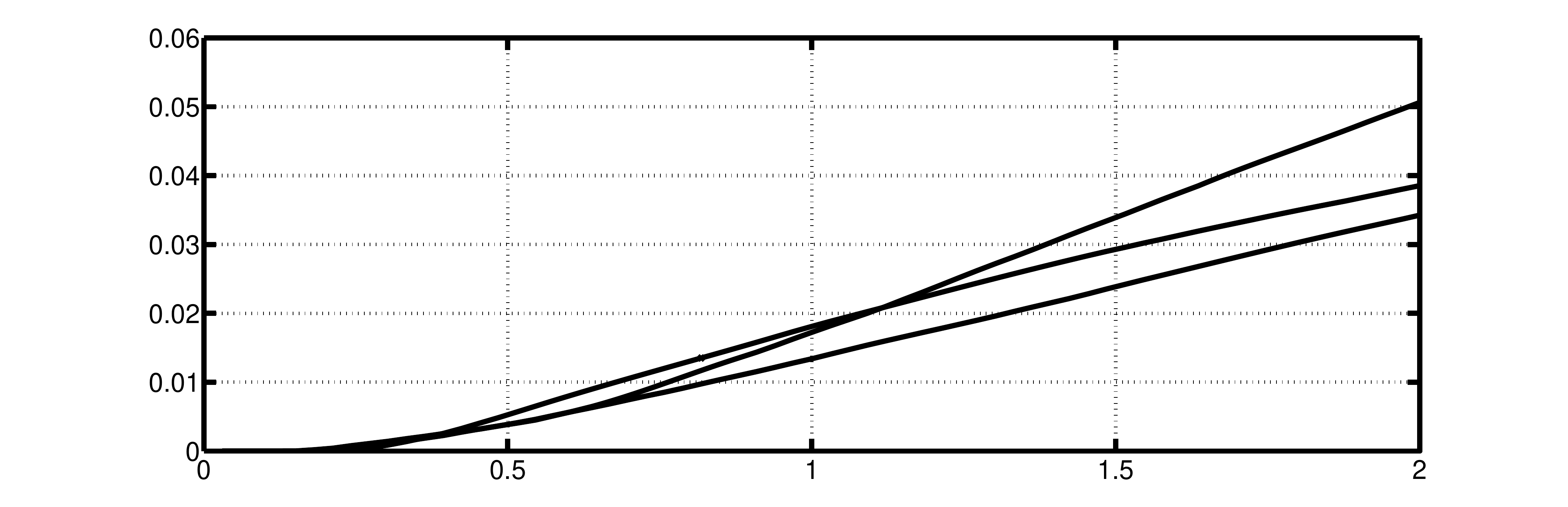}\label{fig_BFSgap}}
\subfigure[Relative gap $\Delta_{avg}$ for random
trees]{\includegraphics[width=0.80\textwidth]{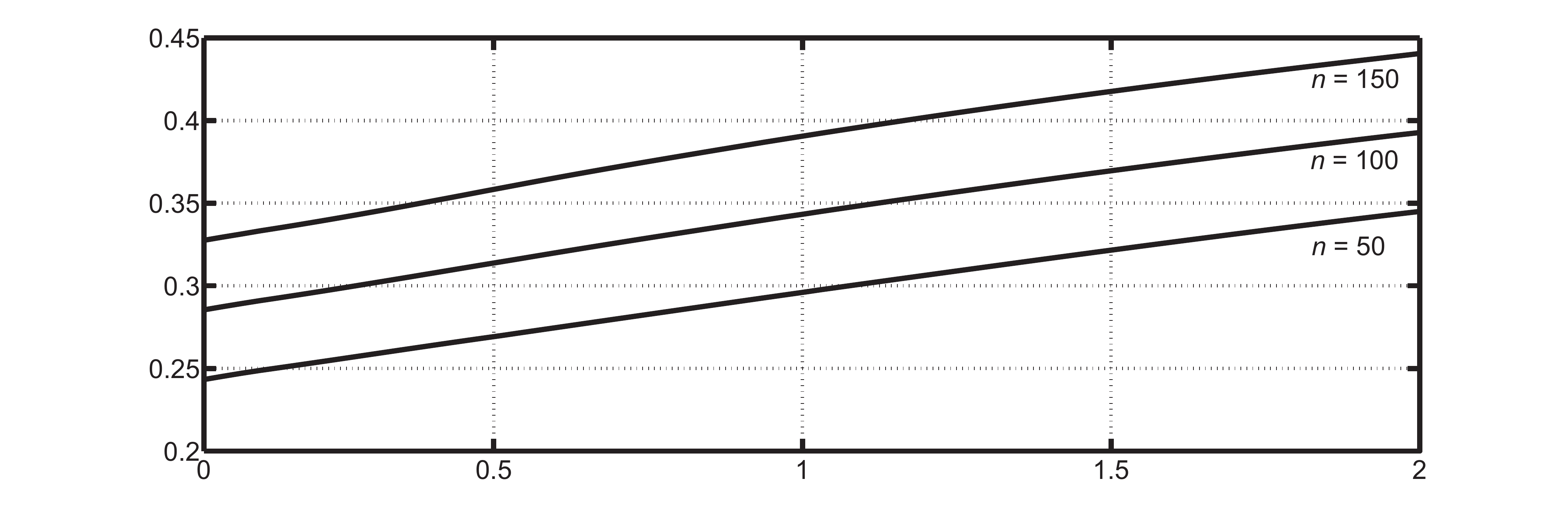}\label{fig_RGas}}
\caption{Gaps \emph{vs} distribution parameter $\beta$: typical relations for
different problem dimension $n$}\label{fig_randomvector}
\end{figure}

From \cite{burkard1985probabilistic} it is known that for QAP the relative gap
between the best and the worst solution tends to zero when the dimension of the
problem increases. For OCN problem, however, the gap depends on the weight
distribution parameter $\beta$. In Figure \ref{fig_BFSgap} three typical relations
are shown between the $BFS$ gap and the weight distribution parameter $\beta$ for
different problem dimension $n$. For $\beta$ being close to zero most weights are
close to unity, and $BFS$ tree is optimal. However, for larger $\beta$ $BFS$ tree is
almost always suboptimal irrespective of the problem dimension. Therefore, even for
random flow matrices the solution of OCN problem can be non-trivial.

The curves in Figure \ref{fig_RGas} show how much we lose in average from choosing a
random tree instead of seeking for a ``good'' tree for $\beta\in[0,2]$. Three
typical curves for different problem dimension $n$ show the significant gap, which
increases when the problem dimension grows.

\subsection{Random flow matrices}
Then the lower bound was tested against a collection of random flow matrices. Again,
100 degree sequences were generated for trees of order from 50 to 250 with degrees
of internal vertices uniformly distributed from 2 to 5. For every degree sequence
$d$ of dimension $n$ a random flow matrix $A=(\mu_{ij})_{i,j=1}^n$ was generated
such that $\mu_{ij}=\mu_{ji}=\rnd^\beta$, where $\rnd$ is a random number uniformly
distributed on $[0,1]$, and $\beta>0$ is a diversity factor. Then each matrix was
loosed to the desired density degree $\sigma$. The lower bound for the best
parameter values, two heuristics, and BFS tree were calculated. As before, the
``best found solution'' was defined as $H^*=\Argmin_{H\in\{H_1,H_2,BFS\}}C_A(H)$,
and the gap of the lower bound was evaluated as $\Delta_{LB}=
\frac{C_A(H^*)-LB_A(d)}{LB_A(d)}$.

\begin{figure}
\subfigure[LB gap $\Delta_{LB}$ \emph{vs} dimension $n$ for weight diversity
$\beta=0,1,2$]{\includegraphics[width=0.5\textwidth]{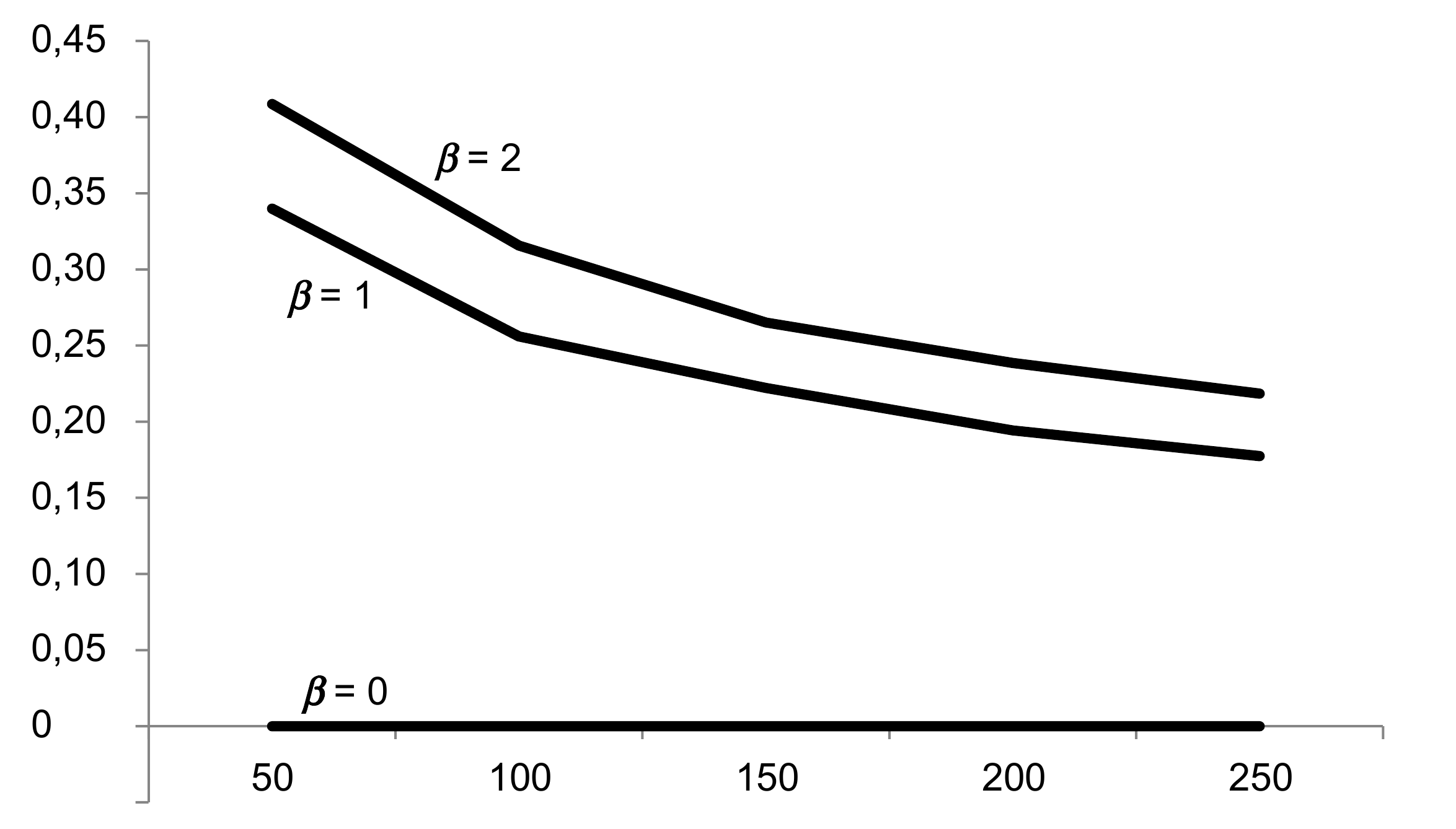}\label{fig_LBgapDim}}
\subfigure[LB gap $\Delta_{LB}$ \emph{vs} flow matrix density $\sigma$ for
$\beta=0,1,2$]{\includegraphics[width=0.50\textwidth]{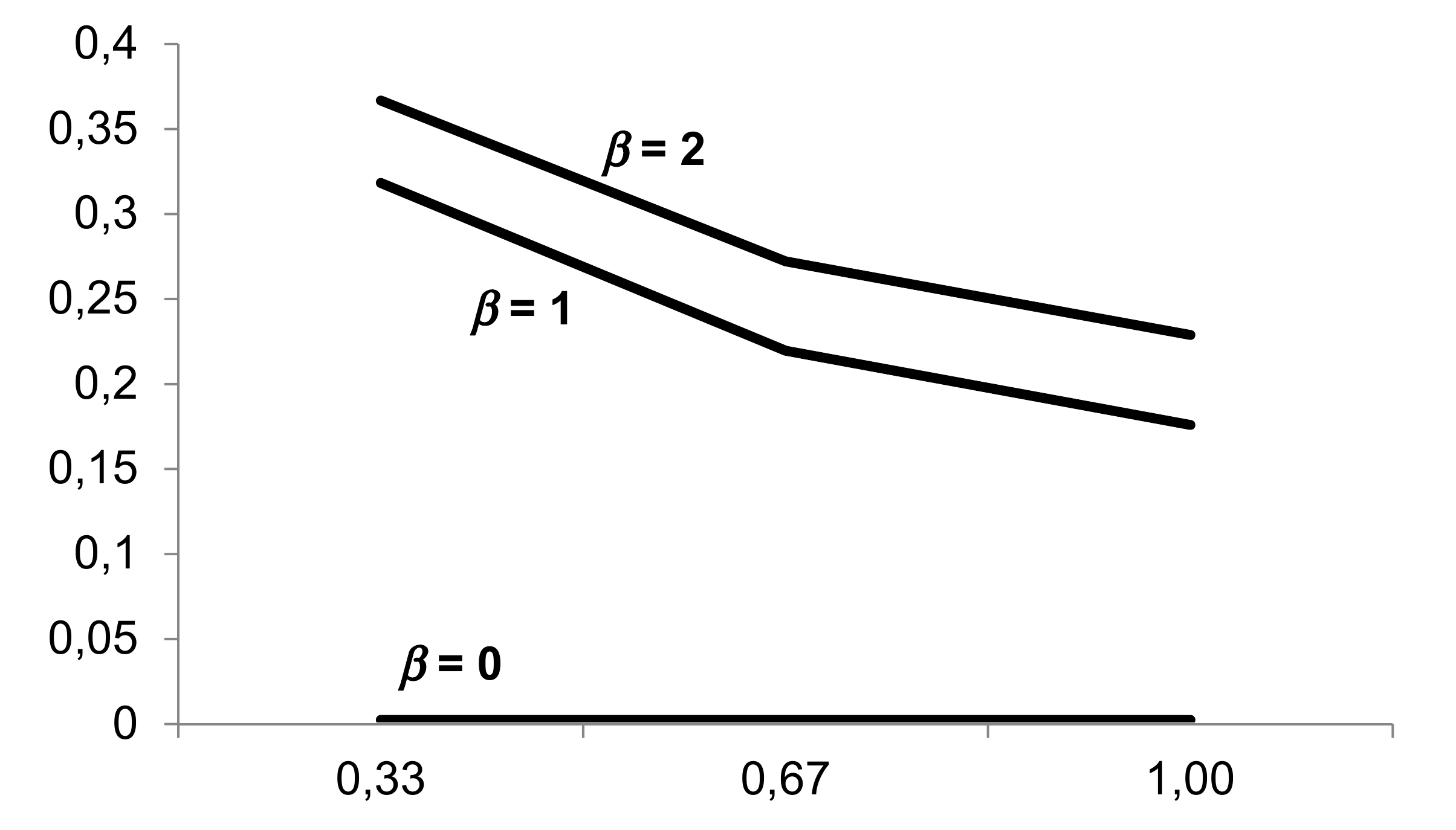}\label{fig_LBgapSpar}}
\subfigure[Random graph gap $\Delta_{avg}$ \emph{vs} problem dimension
$n$]{\includegraphics[width=0.5\textwidth]{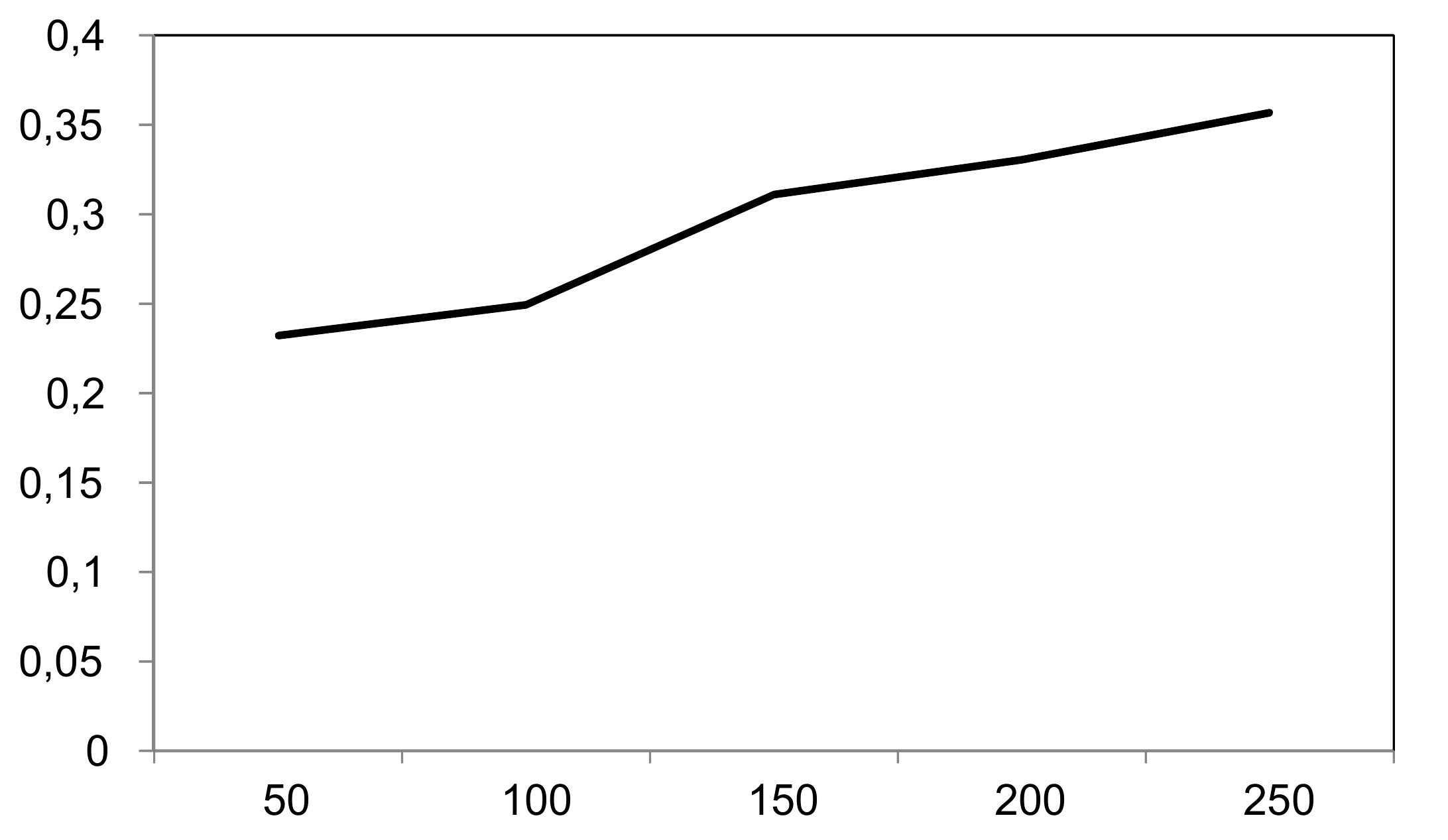}\label{fig_RgapDim}}
\subfigure[BSF tree gap $\Delta_{BFS}$ \emph{vs} problem dimension $n$ for
$\beta=0,1,2$]{\includegraphics[width=0.5\textwidth]{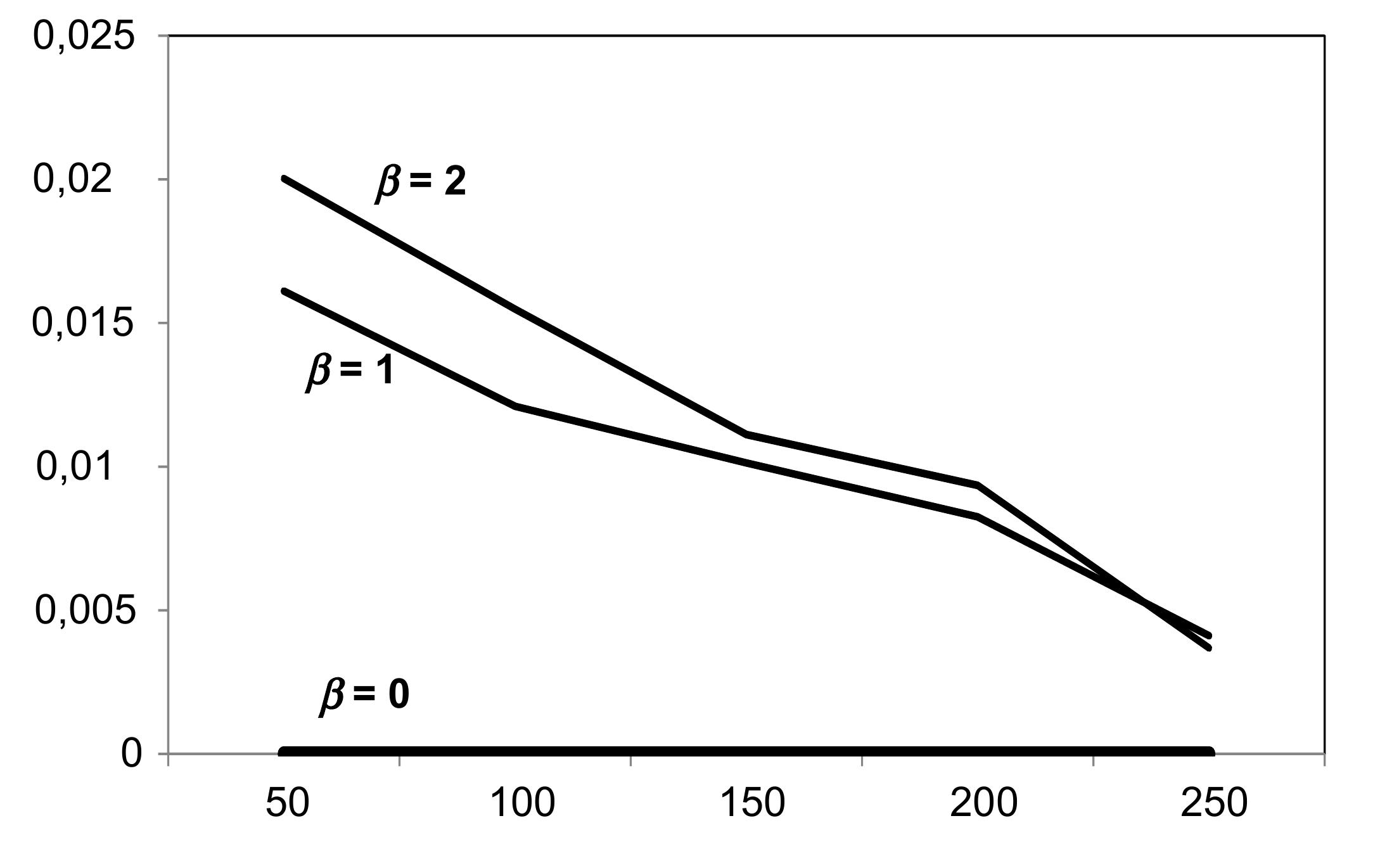}\label{fig_BFSgapDim}}
\subfigure[\# of \textproc{SolveRelaxed} calls \emph{vs} dimension $n$ for
$\beta=0,1,2$
]{\includegraphics[width=0.50\textwidth]{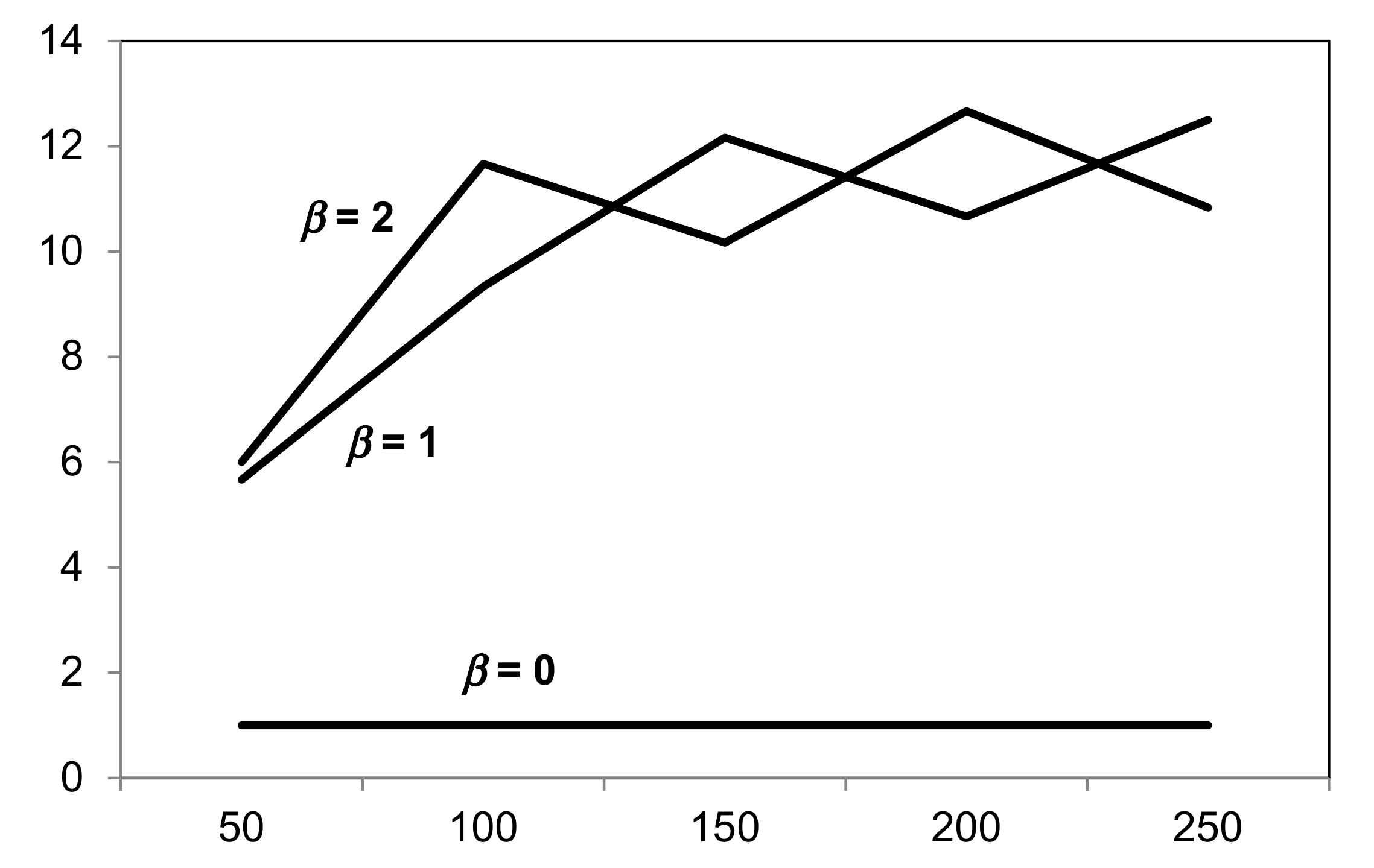}\label{fig_BasisDim}}
\subfigure[Calculation time \emph{vs} problem dimension
$n$]{\includegraphics[width=0.50\textwidth]{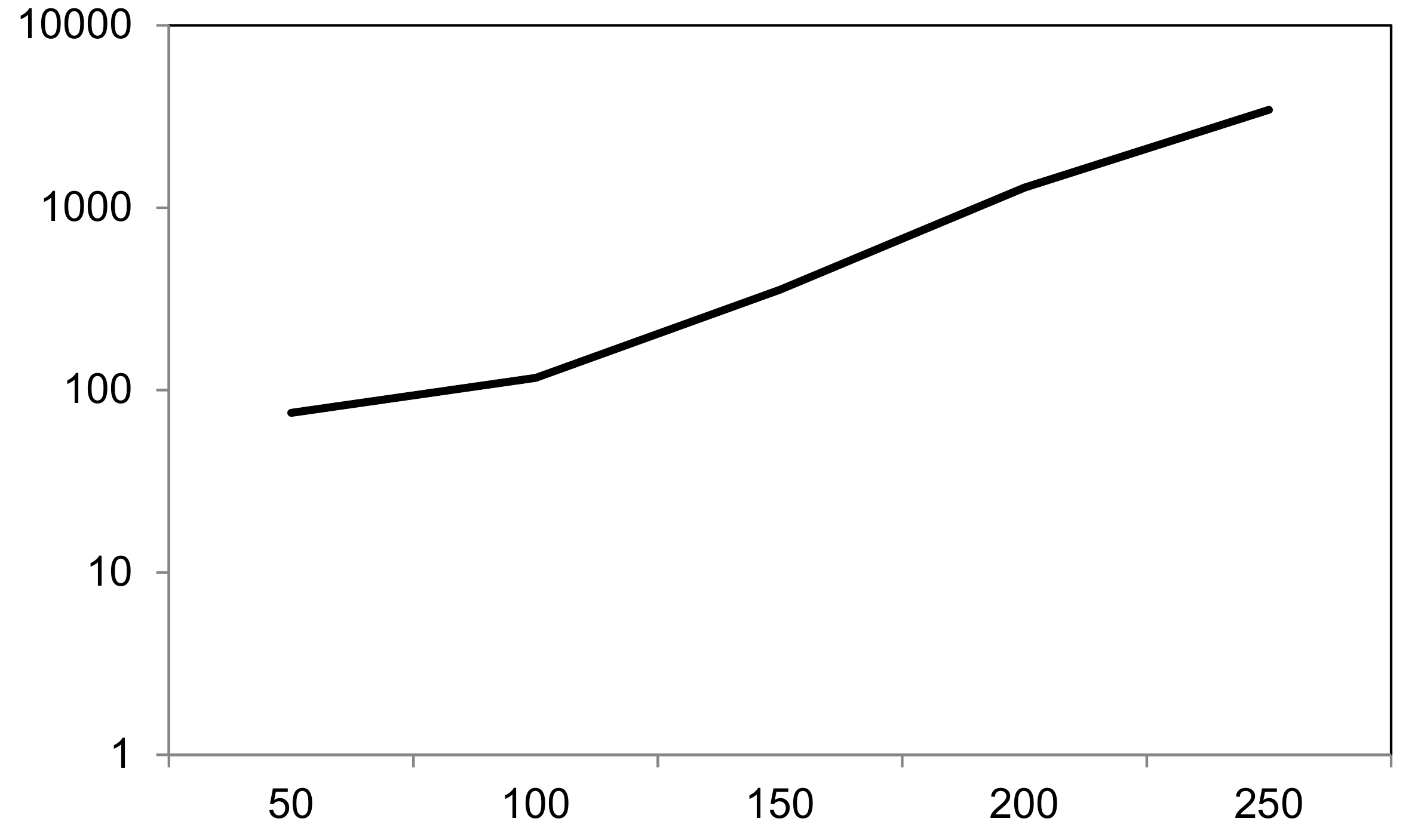}\label{fig_CalcTime}}
\caption{Results of numeric tests for random flow matrices}\label{fig_randommatrix}
\end{figure}

The results are presented in Figure \ref{fig_randommatrix}. For $\beta=0$ we have $A
= J-I$, and OCN problem reduces to the Wiener index minimization whose solution is
BFS tree \cite{wang2008extremal,zhang2008wiener}. The algorithm easily finds optimal
weights (being equal to unity), and the lower-bound gap is equal to zero (see curves
for $\beta=0$ in Figures \ref{fig_LBgapDim} and \ref{fig_LBgapSpar}).

Figure \ref{fig_LBgapDim} shows that the relative gap $\Delta_{LB}$ decreases (and,
hence, the lower bound quality increases) with problem dimension. This effect is
probably due to the random nature of the underlying flow matrices: in a large matrix
the effect of an individual flow is easier to conceal. Again, the average tree gap
$\Delta_{avg}$ in Figure \ref{fig_RgapDim} increases with problem dimension, so the
potential gain from solving OCN problem increases for large-scale problems. At the
same time, Figure \ref{fig_BFSgapDim} shows that BFS tree, the ``universal''
solution that does not depend on the flow matrix, can be very attractive (at least,
when compared to the existing heuristics).

In general, the more diverse are the flows, the lower is the quality of the lower
bound. For example, the sparser matrix enjoys the larger weight diversity, and the
relative gap $\Delta_{LB}$ decreases in matrix density (see Figure
\ref{fig_LBgapSpar}). From Figure \ref{fig_BFSgapDim} we see that the quality of BFS
tree also decreases in weights' diversity $\beta$, and the gain from accounting for
the specific flows' pattern increases.

It is important to note that the average number of calls of \textproc{SolveRelaxed}
function does not increase in problem dimension (see Figure \ref{fig_BasisDim}), so
even for large-scale problems wee need not consider bulky constraint sets in
(\ref{eq_sup_huffman}). At the same time, several large-scale convex problems are
solved inside \textproc{SolveRelaxed}, which requires more calculus when the problem
dimension grows. The computation time increases rapidly in problem dimension $n$
(see Figure \ref{fig_CalcTime}). However, problems with several hundreds terminals
are still solved in reasonable time (see details in Figure \ref{fig_CalcTime}) on a
laptop (we used Lenovo\texttrademark Thinkpad\textsuperscript{\copyright} with
Intel\texttrademark Core i5 2.3GHz).

\subsection{Real-world datasets}
    Different free data sources from transportation industry (public transport and airline statistics reports) and
    demography (migration reports) were used to build several
    real-world flow matrices with various size and flows' pattern. Below we briefly characterize all sources. Information about all elicited datasets is consolidated in Table
    \ref{tab_dataset}. In all cases, we symmetrize obtained origin-destination (OD) matrices to obtain a symmetric flow matrix.
    The flow matrices can be downloaded from \url{http://www.mtas.ru/upload/ODmatrices.zip}.

    \begin{enumerate}\sloppy
    \item \textbf{London Tube and Rail Transport (LTRT)}\\
    It is possible to travel on Tube, DLR, London Overground, TfL Rail and most National Rail services using
    contactless or Oyster card to pay. Two data sets located at \url{https://tfl.gov.uk/maps/track/dlr} provide information about the traffic between London Tube
    and Rail Stations based upon the card touch-in/touch-out information.

    \item \textbf{Queensland Government Data --TransLink OD trips (TransLink)}\\
    Several datasets were derived from the Queensland State's Government,
    Australia,
    \url{https://data.qld.gov.au/dataset/go-card-transaction-data/resource/8a99a319-6870-4945-b87e-e58b178deae3},
    storing data about trip count for many transportation modes and
    carriers.

    \item\textbf{Greater Cambridge ANPR Data: OD Reports (ANPR)}\\
    These origin-to-destination reports are derived from the Automatic Number Plate Recognition (ANPR)
    camera traffic survey undertaken in 2017 across the Cambridge area from June 10 to 17.
    The reports provide information on the first and the last cameras
triggered on vehicle journeys across the city. We summarize the data into an OD
matrix.

    \item \textbf{The Air Carrier Statistics database -- T-100 Segment (T100)}\\
The source located at \url{https://www.transtats.bts.gov/Fields.asp?Table_ID=293}
contain domestic and international T-100 segment data reported by U.S. and foreign
air carriers and non-stop segment data by aircraft type and service class for
transported passengers, freight and mail, available capacity, scheduled departures,
departures performed, aircraft hours, and load factor. Flights with both origin and
destination in a foreign country are not included. OD matrix is built using the
fields ``OriginAirportID'', ``DestAirportID'', ``Passengers'', and
``UniqueCarrier''.

    \item \textbf{Airline Origin and Destination Survey (US Air)}\\
    The Airline Origin and Destination Survey (\url{https://data.world/us-dot-gov/02210b59-4330-440d-acf4-d4fb276f1d74}) is a 10\% sample of
    airline tickets from reporting carriers collected by the U.S. Office of
    Airline Information of the Bureau of Transportation Statistics in the first quarter of 1993.
    Data includes origin, destination and other itinerary details
    of passengers transported. This database is used to determine
    air traffic patterns, air carrier market shares and passenger flows.
    We analyze only fields ``OriginAirportID'', ``DestAirportID'', ``Coupons''.
    If it were several airports in an itinerary, we take the first airport as origin an the last as destination.
    If the first airport coincide with the last, we split the itinerary on two itineraries:
    from the first airport to the penultimate and from the penultimate to the last one.

    \item\textbf{Canada Aircraft Movement Statistics (Canadian)}

    The survey located at \url{http://www23.statcan.gc.ca/imdb/p2SV.pl?Function=getSurvey&SDDS=2715} provides
    estimates of aircraft movements in Canada. The source table
    contains the hyphen-separated pair of cities and the passenger flow between these cities.

    \item \textbf{EU Country to Country Migration (EU Migration)}\\
    \url{https://www.imi.ox.ac.uk/data/demig-data/demig-c2c-data}

    The DEMIG C2C (country-to-country) database contains bilateral migration flow data for 34
    reporting countries and from up to 236 countries over the 1946--2011 period.
    It includes data for inflows, outflows and net flows, respectively for citizens,
    foreigners and/or citizens and foreigners combined, depending on the reporting countries.
    We take ``Reporting country'', ``Countries'', and ``Value'' columns for both
genders.

    \item \textbf{U.S. Census Bureau Migration Reports (US Migration)}\\
    The U.S. Census Bureau has been releasing county-to-county and county/minor civil division (MCD)-to-county/MCD migration flow
    estimates based on the American Community Survey (ACS) since 2012. We use the columns $i=$``FIPS County Code of Geography A'', $j=$``FIPS County Code of Geography B'',
    $f=$``Flow from Geography B to Geography A'', ``Counterflow from Geography A to Geography B1'' to construct the
    symmetric matrix $c_{ij} = (f_{ij}+f_{ji})/2$.
\end{enumerate}
\begin{table}\sloppy\centering
\caption{Data sets used to build real-world flow matrices}
\label{tab_dataset}\begin{small}\tabcolsep=0.11cm\begin{tabular}{|llcc|}
  \hline
  Source & Dataset & Abbrev. & Dimension \\
  \hline
  LTRT & The London Underground Limited operator  & LUL & 266 \\     
  LTRT & The Docklands Light Railway light metro system & DLR & 61 \\ 
  TransLink & All carriers in June 2017 & TL & 723 \\
  TransLink & One week of June 2017 for the carrier ``BCC Ferries'' & BCC & 20 \\
  TransLink & One week of June 2017 for the carrier ``Sunbus'' & Sunbus & 846 \\
  TransLink & Carrier ``Park Ridge Transit'' in June 2017& PRT & 498 \\
  TransLink & Carrier ``Mt Gravatt Bus Service'' in June 2017& MGBS & 364 \\
  TransLink & Carrier ``Queensland Rail'' in June 2017& QR & 154 \\
  ANPR  & Summary data for June 10, 2017 & ANPR & 91 \\
  T100 & Carrier ``Hawaiian Airlines Inc'' in January 2017& HA & 29 \\
  T100 & Carrier ``Compass Airlines'' in January 2017& CA & 55 \\
  US Air & Carrier ``America West Airlines Inc.'' (IATA code HP) & HP & 105 \\
  US Air & Carrier ``Trans World Airways LLC'' (IATA code TW) & TW & 176 \\
  US Air & Carrier ``US Airways Inc.'' (IATA code US) & US & 269 \\
  US Air & Carrier ``Midwest Express Airlines'' (IATA code YX) & YX & 59 \\
  Canadian & The annual report for 2005 & Canadian & 72 \\
  EU Migration & EU to EU migration in 2007 by country & EU & 13 \\
  US Migration & Migration between counties of Alabama in 2014 & Alabama & 67 \\
  \hline
\end{tabular}
\end{small}
\end{table}

\begin{figure}
\subfigure[ANPR (public
transport)]{\includegraphics[width=0.50\textwidth]{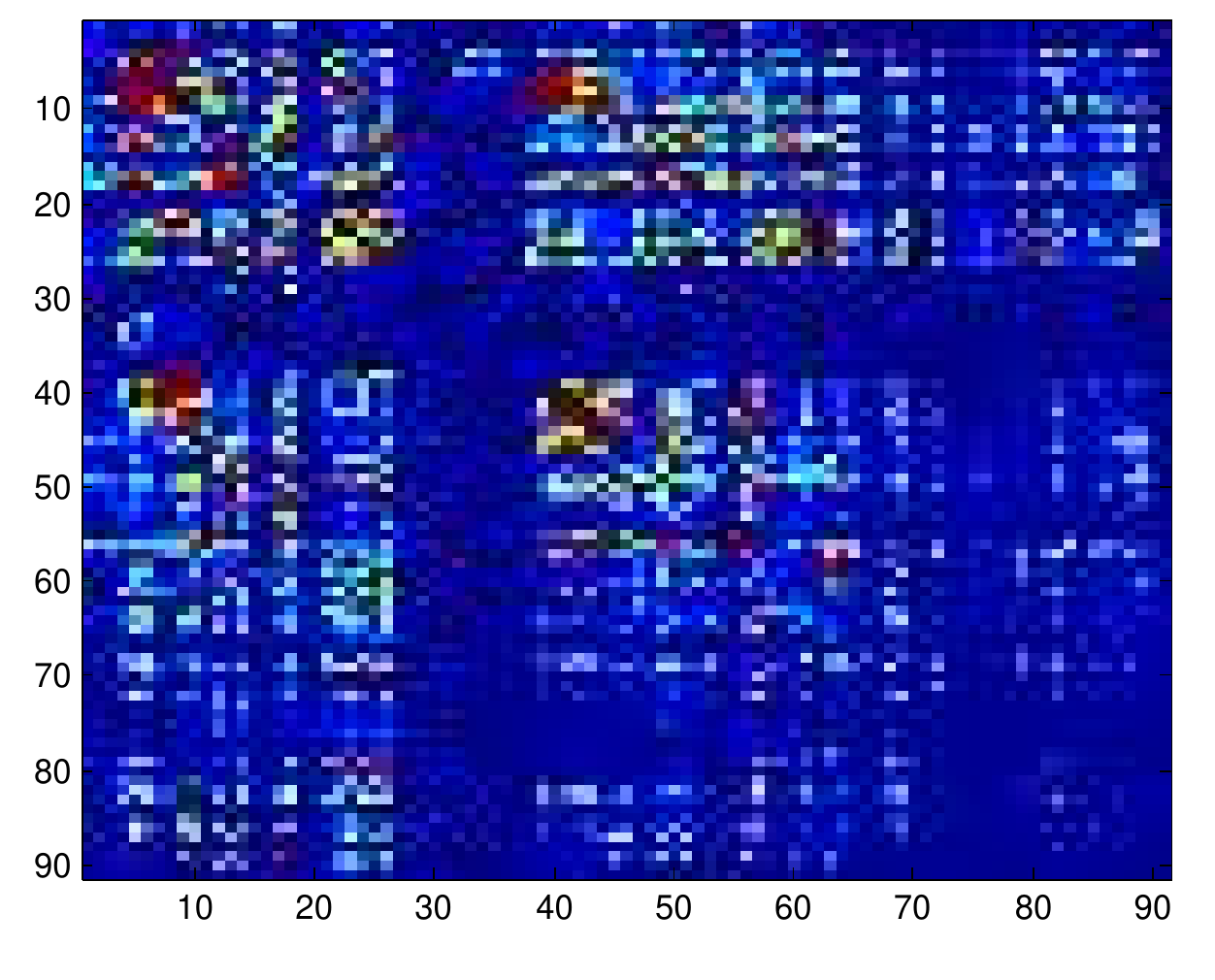}\label{fig_ANPR}}
\subfigure[QR (public
transport)]{\includegraphics[width=0.50\textwidth]{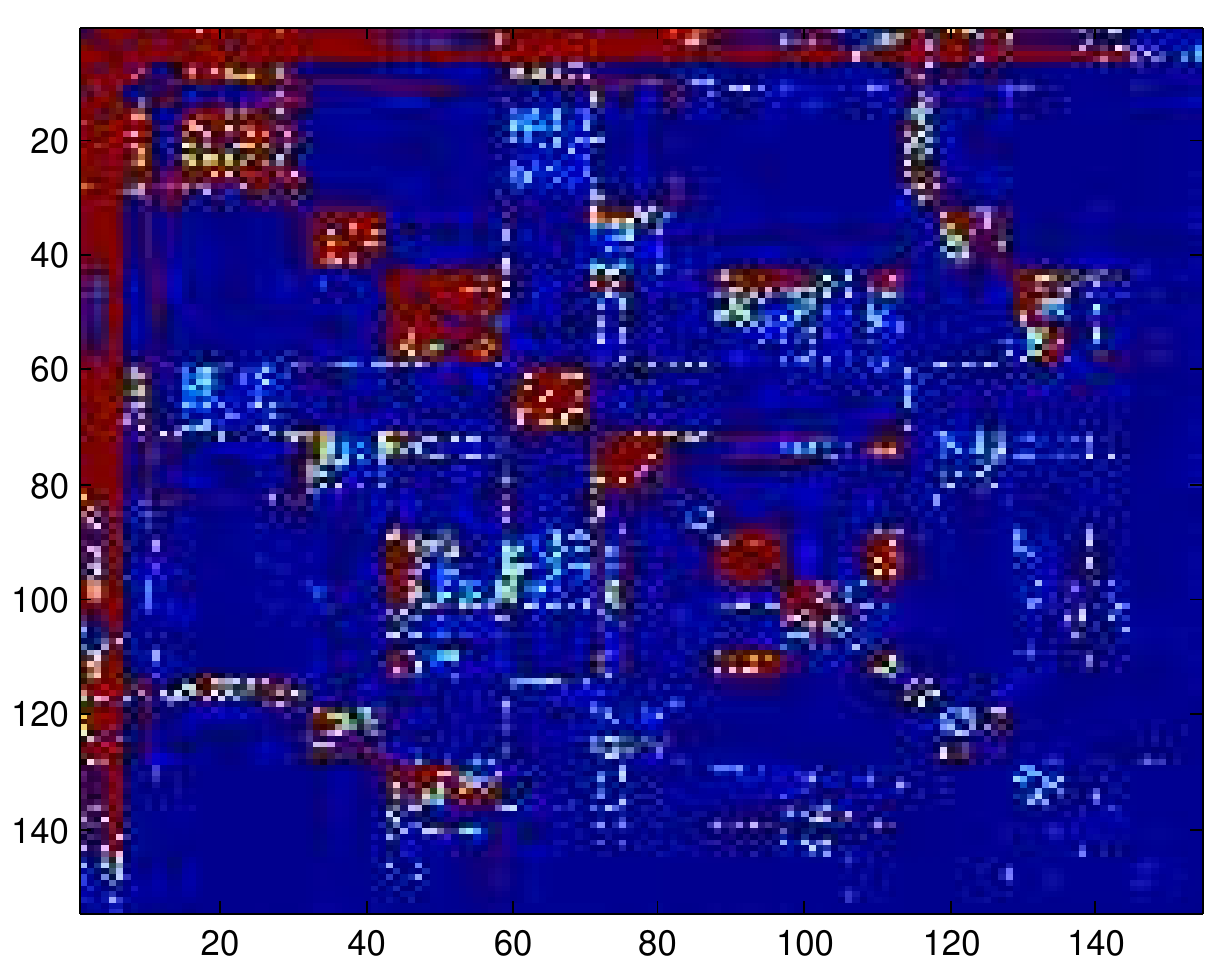}\label{fig_QR}}\\
\subfigure[Alabama
(migration)]{\includegraphics[width=0.50\textwidth]{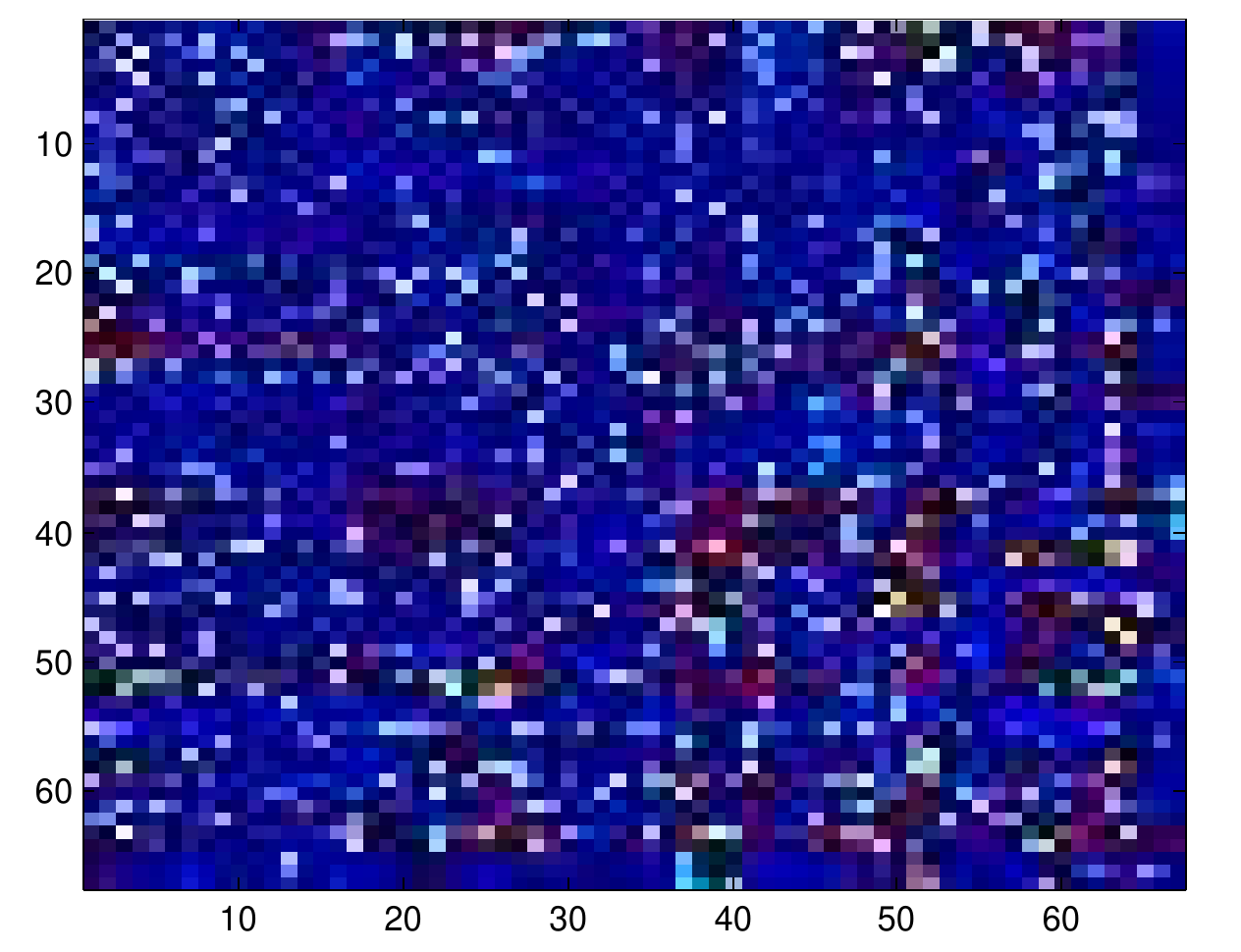}\label{fig_Alabama}}
\subfigure[HP
(aviation)]{\includegraphics[width=0.50\textwidth]{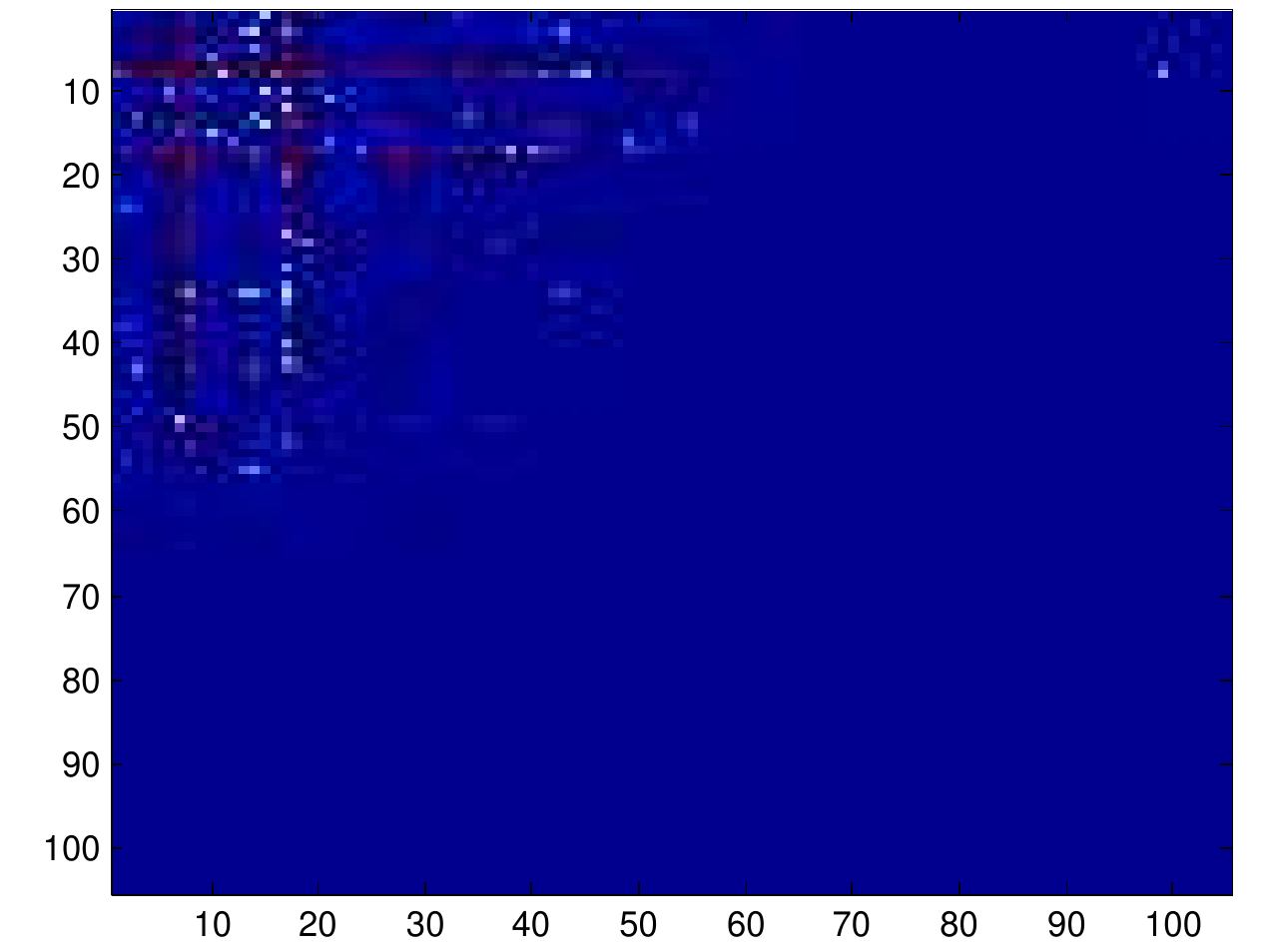}\label{fig_HP}}
\caption{Typical OD matrices before symmetrization (flow intensity grows from blue
to red)}\label{fig_datasets}
\end{figure}

Several typical flow patterns are presented in Figure \ref{fig_datasets}.
Transportation and migration datasets were used because of their availability,
although we clearly understand that minimizing the number-of-edges graph distance
over the set of trees is not of much practical interest for them.

\begin{figure}
{\includegraphics[width=0.99\textwidth]{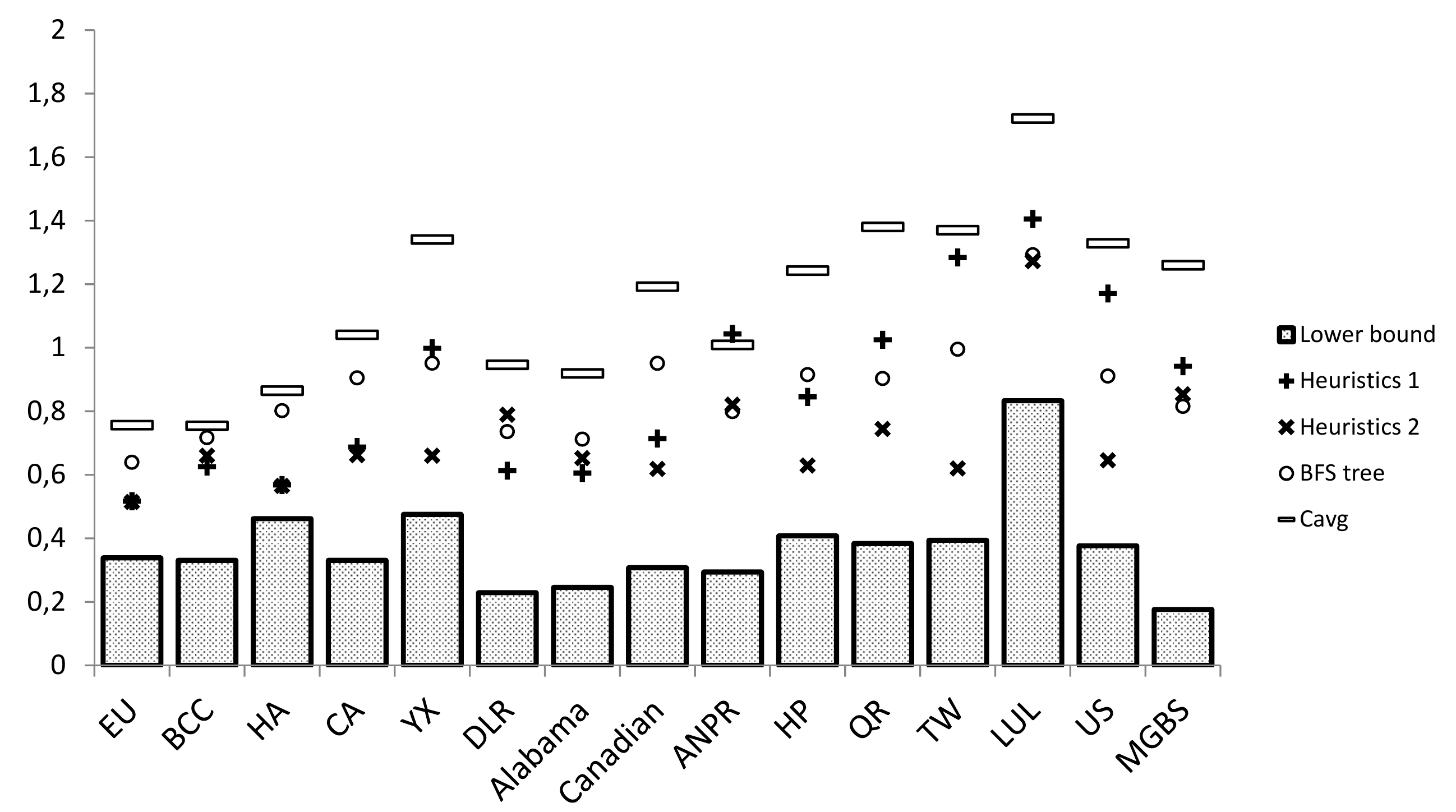}} \caption{Results for real
datasets (normalized to the sizing factor $n^2\ln n$). The lower bound is depicted
with the height of a dotted rectangle, ``$+$'' shows the cost of
\textproc{Heuristics1} tree, while ``$\times$'' stands for \textproc{Heuristics2}
tree. Circle ``$\circ$'' shows the cost of BFS tree, and the horizontal bar points
out the cost of the average tree.}\label{fig_realdata}
\end{figure}

For each flow matrix the degree sequence was generated with degrees of internal
vertices uniformly distributed from 2 to 5. The lower bound, two upper bounds, BFS
tree, and the average tree cost $C_{avg}$ were calculated. The results are presented
in Table \ref{tab_results} and in Figure \ref{fig_realdata}.

\begin{table}[htbp]
\sloppy \centering \caption{Results for real-world datasets. The cost of the best
found tree is marked with bold}
\label{tab_results}
\tabcolsep=0.11cm\begin{small}
    \begin{tabular}{|lc||cccccccc|}
    \toprule
    \multicolumn{1}{|c}{Dataset} & Dimension & \multicolumn{1}{c}{$\Delta_{LB}$} & \multicolumn{1}{c}{$\Delta_{BFS}$} & \multicolumn{1}{c||}{$\Delta_{avg}$} & \multicolumn{1}{c}{LB} & \multicolumn{1}{c}{Heur. 1} & \multicolumn{1}{c}{Heur. 2} & \multicolumn{1}{c}{BFS tree} & \textit{$C_{avg}$} \\
    \midrule
    EU    & 13    & \multicolumn{1}{c}{52\%} & \multicolumn{1}{c}{24\%} & \multicolumn{1}{c||}{47\%} & \multicolumn{1}{r}{0.339} & \multicolumn{1}{r}{0.517} & \multicolumn{1}{r}{\textbf{0.515}} & \multicolumn{1}{r}{0.639} & \multicolumn{1}{r|}{0.756} \\
    BCC   & 20    & \multicolumn{1}{c}{89\%} & \multicolumn{1}{c}{15\%} & \multicolumn{1}{c||}{20\%} & \multicolumn{1}{r}{0.331} & \multicolumn{1}{r}{\textbf{0.626}} & \multicolumn{1}{r}{0.660} & \multicolumn{1}{r}{0.717} & \multicolumn{1}{r|}{0.754} \\
    HA    & 29    & \multicolumn{1}{c}{22\%} & \multicolumn{1}{c}{42\%} & \multicolumn{1}{c||}{53\%} & \multicolumn{1}{r}{0.462} & \multicolumn{1}{r}{0.568} & \multicolumn{1}{r}{\textbf{0.565}} & \multicolumn{1}{r}{0.802} & \multicolumn{1}{r|}{0.865} \\
    CA    & 55    & \multicolumn{1}{c}{100\%} & \multicolumn{1}{c}{37\%} & \multicolumn{1}{c||}{57\%} & \multicolumn{1}{r}{0.331} & \multicolumn{1}{r}{0.687} & \multicolumn{1}{r}{\textbf{0.662}} & \multicolumn{1}{r}{0.905} & \multicolumn{1}{r|}{1.041} \\
    YX    & 59    & \multicolumn{1}{c}{39\%} & \multicolumn{1}{c}{44\%} & \multicolumn{1}{c||}{103\%} & \multicolumn{1}{r}{0.476} & \multicolumn{1}{r}{0.998} & \multicolumn{1}{r}{\textbf{0.660}} & \multicolumn{1}{r}{0.951} & \multicolumn{1}{r|}{1.340} \\
    DLR   & 61    & \multicolumn{1}{c}{168\%} & \multicolumn{1}{c}{20\%} & \multicolumn{1}{c||}{54\%} & \multicolumn{1}{r}{0.229} & \multicolumn{1}{r}{\textbf{0.613}} & \multicolumn{1}{r}{0.789} & \multicolumn{1}{r}{0.736} & \multicolumn{1}{r|}{0.946} \\
    Alabama & 67    & \multicolumn{1}{c}{147\%} & \multicolumn{1}{c}{18\%} & \multicolumn{1}{c||}{52\%} & \multicolumn{1}{r}{0.245} & \multicolumn{1}{r}{\textbf{0.605}} & \multicolumn{1}{r}{0.652} & \multicolumn{1}{r}{0.712} & \multicolumn{1}{r|}{0.919} \\
    Canadian & 72    & \multicolumn{1}{c}{101\%} & \multicolumn{1}{c}{54\%} & \multicolumn{1}{c||}{93\%} & \multicolumn{1}{r}{0.307} & \multicolumn{1}{r}{0.714} & \multicolumn{1}{r}{\textbf{0.619}} & \multicolumn{1}{r}{0.951} & \multicolumn{1}{r|}{1.192} \\
    ANPR  & 91    & \multicolumn{1}{c}{171\%} & \multicolumn{1}{c}{0\%} & \multicolumn{1}{c||}{26\%} & \multicolumn{1}{r}{0.294} & \multicolumn{1}{r}{1.044} & \multicolumn{1}{r}{0.821} & \multicolumn{1}{r}{\textbf{0.799}} & \multicolumn{1}{r|}{1.009} \\
    HP    & 105   & \multicolumn{1}{c}{54\%} & \multicolumn{1}{c}{45\%} & \multicolumn{1}{c||}{97\%} & \multicolumn{1}{r}{0.408} & \multicolumn{1}{r}{0.845} & \multicolumn{1}{r}{\textbf{0.629}} & \multicolumn{1}{r}{0.915} & \multicolumn{1}{r|}{1.243} \\
    QR    & 154   & \multicolumn{1}{c}{94\%} & \multicolumn{1}{c}{21\%} & \multicolumn{1}{c||}{85\%} & \multicolumn{1}{r}{0.384} & \multicolumn{1}{r}{1.025} & \multicolumn{1}{r}{\textbf{0.745}} & \multicolumn{1}{r}{0.903} & \multicolumn{1}{r|}{1.380} \\
    TW    & 176   & \multicolumn{1}{c}{57\%} & \multicolumn{1}{c}{61\%} & \multicolumn{1}{c||}{121\%} & \multicolumn{1}{r}{0.394} & \multicolumn{1}{r}{1.284} & \multicolumn{1}{r}{\textbf{0.620}} & \multicolumn{1}{r}{0.995} & \multicolumn{1}{r|}{1.370} \\
    LUL   & 266   & \multicolumn{1}{c}{53\%} & \multicolumn{1}{c}{2\%} & \multicolumn{1}{c||}{35\%} & \multicolumn{1}{r}{0.833} & \multicolumn{1}{r}{1.405} & \multicolumn{1}{r}{\textbf{1.272}} & \multicolumn{1}{r}{1.293} & \multicolumn{1}{r|}{1.721} \\
    US    & 269   & \multicolumn{1}{c}{71\%} & \multicolumn{1}{c}{41\%} & \multicolumn{1}{c||}{106\%} & \multicolumn{1}{r}{0.377} & \multicolumn{1}{r}{1.171} & \multicolumn{1}{r}{\textbf{0.646}} & \multicolumn{1}{r}{0.911} & \multicolumn{1}{r|}{1.329} \\
    MGBS  & 364   & \multicolumn{1}{c}{363\%} & \multicolumn{1}{c}{0\%} & \multicolumn{1}{c||}{54\%} & \multicolumn{1}{r}{0.176} & \multicolumn{1}{r}{0.942} & \multicolumn{1}{r}{0.853} & \multicolumn{1}{r}{\textbf{0.815}} & {1.259} \\
    PRT   & 498   & \multicolumn{8}{c|}{Stopped after several hours of computation} \\
    TL    & 723   & \multicolumn{8}{c|}{Stopped after several hours of computation} \\
    Sunbus & 846   & \multicolumn{8}{c|}{Stopped after several hours of computation} \\
    \bottomrule
    \end{tabular}
\end{small}
\end{table}

Table \ref{tab_results} shows that, distinct to random flows (see Figure
\ref{fig_BFSgapDim}), BFS tree can have unacceptable quality for real datasets.
Also, in most cases, \textproc{Heuristics2} gives the best tree. Therefore, in spite
of the complexity of \textproc{MAximizeLB} procedure, it provides the highly
valuable information for heuristic algorithm construction.

From Figure \ref{fig_realdata} wee see that the quality of the lower bound (the
value of the relative gap $\Delta_{LB}$) may vary in a wide range: from the modest
gap $\Delta_{LB}=22\%$ for ``HA'' dataset to the huge gap $\Delta_{LB}=363\%$ for
``MGBS'' dataset. In the latter case the lower bound becomes almost uninformative
(although it is still twice as big as a trivial lower bound, the sum of all flows).
At the same time, we do not see the quality of the lower bound to decrease with the
problem dimension. So, good quality of the lower bound can be expected for bigger
samples, at least for some application areas. Potentially, after the careful
optimization of the algorithm, the best parameters of the lower bound can be
calculated for OCN problems with thousand terminals or more.

%

\section{Conclusion}

An optimal connecting network (OCN) has the minimum possible weighted sum of
distances between pairs of its vertices among all admissible networks. Weights of
vertex pairs are given by a flow matrix $A$. In general, finding OCN is a complex
problem of combinatory optimization.

In this article a lower bound estimate is constructed for the cost of an optimal
connecting tree with the given degree sequence. The lower bound is parameterized by
two vectors, $\alpha\in \mathbb{R}^n$ and $\mu\in \mathbb{R}^n_+$. The problem of
finding the best combination of parameter values reduces to the non-convex
semidefinite problem, for which an algorithm is proposed. The algorithm solves the
non-convex problem through a series of its convex relaxations.

Although the optimization problem involves a (rather demanding) semidefinite
constraint and several quadratic constraints with dense matrices, numeric tests show
that the lower bound can be calculated in reasonable time (minutes on a PC) for
trees with several hundreds vertices. However, calculation of the lower bound for
huge trees with thousands vertices is still an open problem, which can be the
subject of future research. At the same time, if we do not insist on the best
parameter values and are satisfied with any admissible $\alpha$ and $\mu$,
calculation time can be considerably decreased by increasing tolerance parameter
$\delta$ in Listing \ref{lst:alg}.

The quality of the lower bound depends on how accurately flow matrix $A$ can be
approximated by the sum of the diagonal matrix $\diag(\alpha)$ and the non-negative
rank-one matrix $\mu\mu^\top$. It is shown in Section \ref{sec_numeric} that for $A$
having rank one we have the perfect approximation, and the lower bound is equal to
the optimal tree cost. In this case every terminal can be endowed with non-negative
weight $\mu_i$, $i=1,...,n$, and the flow between the terminals $i$ and $j$ is
written as $\mu_{ij} = \mu_i\mu_j$. Weights of terminals are explained by the
following simplistic model. Let us assume that the $i$-th terminal is \emph{active}
at a given period of time with probability proportional to its weight $\mu_i$. If
active, a terminal sends a unique piece of information to all terminals being active
at this moment. If all terminals are independent, then the average volume of
information circulating between terminals $i$ and $j$ is proportional to
$\mu_i\mu_j$.

Many real-world flow patterns, however, are far from this model, and the lower bound
may sometimes have poor quality. It is an open question, which flow matrix is the
least convenient for approximation by a rank-one matrix, and, hence, for which flow
matrix the lower bound has the least quality. These results may be used when
developing the new lower bounds with the better guaranteed quality.


The strategic direction of research, however, is connected with generalizing the
approach to the general networks with loops.

\section*{Funding}

This work was supported by the Russian Foundation for Basic Research (RFBR)
[18-07-01240].

\bibliographystyle{plain}
\bibliography{GoubkoKuznetsov}

\appendix

\section{Properties of bilinear matrix inequality $xx^\top-A \succeq 0$}\label{sec_appendix}

In this appendix properties are studied of the set
$$X_A := \left\{x\in\mathbb{R}^n: xx^\top-A \succeq 0\right\}$$
where $A$ is a symmetric real $n\times n$ matrix.

Recall that with $\lambda_i(A), i=1,...,n$ we denote (real) eigenvalues of real
symmetric matrix $A$ listed in the descending order while $u^{(i)}(A)$ standing for
the corresponding eigenvectors. Let  $\mathbb{S}^{n-1} := \left\{x\in \mathbb{R}^n :
x^\top x=1\right\}$ denote the unit sphere in $\mathbb{R}^n$.

\begin{lemma}\label{lemma_lambda2_empty}
$X_A\neq \emptyset$ if and only if $\lambda_2(A)\leq 0$.

\begin{proof}
Vector ${x}$ belongs to $X_A$ if and only if for any vector ${z}\in
\mathbb{S}^{n-1}$ inequality ${z}^\top \left({xx}^\top-A\right){z}\ge 0$ holds.
Consequently, $X_A = \emptyset$ if and only if for any ${x}\in \mathbb{R}^n$ there
exists such ${z}\in \mathbb{S}^{n-1}$ that ${z}^\top\left(A-{xx}^\top\right){z} >
0$. In the other words, $X_A = \emptyset$ when
\begin{equation}\label{eq_X_is_empty}
\inf_{a\ge 0} \min_{{x}\in \mathbb{S}^{n-1}} \max_{{z}\in \mathbb{S}^{n-1}}
{z}^\top\left(A-a{xx}^\top\right){z} > 0.
\end{equation}

It is clear that the left-hand side of inequality (\ref{eq_X_is_empty}) will not
increase if we narrow the maximization area, and, therefore,
\begin{multline}\label{eq_ge_lambda2}
\inf_{a\ge 0} \min_{{x}\in \mathbb{S}^{n-1}} \max_{{z}\in
\mathbb{S}^{n-1}} {z}^\top\left(A-a\cdot{xx}^\top\right){z} \ge\\
\ge \inf_{a\ge 0} \min_{{x}\in \mathbb{S}^{n-1}} \max_{{z}\in \mathbb{S}^{n-1}, {z}
\bot {x}}{z}^\top\left(A-a{xx}^\top\right){z}=\\
=\inf_{a\ge 0} \min_{{x}\in \mathbb{S}^{n-1}} \max_{{z}\in \mathbb{S}^{n-1}, {z}
\bot {x}} {z}^\top A {z}=\lambda_2(A).
\end{multline}
The last equality follows from the Courant-Fischer theorem , which says that
$$\lambda_2(A) = \min_{{x}\in \mathbb{S}^{n-1}} \max_{{z}\in \mathbb{S}^{n-1}, {z}\bot {x}} {z}^\top A {z}.$$

On the other hand, the left-hand side of inequality (\ref{eq_X_is_empty}) will not
decrease if minimization over $x\in \mathbb{S}^{n-1}$ if replaced with the concrete
$x = u^{(1)}(A)$:
\begin{multline}
\inf_{a\ge 0} \min_{{x}\in \mathbb{S}^{n-1}} \max_{{z}\in
\mathbb{S}^{n-1}} {z}^\top\left(A-a\cdot xx^\top\right){z} \le\\
\le \inf_{a\ge 0} \max_{{z}\in \mathbb{S}^{n-1}}
{z}^\top\left(A-au^{(1)}(A)u^{(1)}(A)^\top\right){z}=\\
= \inf_{a\ge 0} \lambda_1\left(A-au^{(1)}(A)u^{(1)}(A)^\top\right).
\end{multline}
The last equality also follows from Courant-Fischer theorem.

The spectrum of matrix $A-a\cdot u^{(1)}(A)u^{(1)}(A)^\top$ differs from that of
matrix $A$ only in one component: the eigenvalue $\lambda_1(A)$ is replaced with
$\lambda_1(A)-a$, and so,
\begin{equation}\label{eq_le_lambda2}
\inf_{a\ge 0} \lambda_1\left(A-au^{(1)}(A)u^{(1)}(A)^\top\right)=\inf_{a\ge 0}
\max\left[\lambda_1(A)-a, \lambda_2(A)\right]=\lambda_2(A).
\end{equation}

From inequalities (\ref{eq_ge_lambda2}) and (\ref{eq_le_lambda2}) it follows that
the left-hand side of inequality (\ref{eq_X_is_empty}) is equal to $\lambda_2(A)$,
and the inequality $\lambda_2(A)>0$ is necessary and sufficient for inequality
(\ref{eq_X_is_empty}) to be valid, which, in turn, implies that $X_A$ is empty.
\end{proof}
\end{lemma}

\begin{lemma}\label{lemma_negative_semidefinite}
If matrix $A$ is negative definite, then $X_A=\mathbb{R}^n$.

\begin{proof}
The proof follows immediately from positive semidefiniteness of matrix ${xx}^\top$
for arbitrary ${x}\in \mathbb{R}^n$.
\end{proof}
\end{lemma}

\begin{lemma}\label{lemma_perron_vector}
if $X_A$ is not empty, then ${x} := \sqrt{\lambda_1(A)}u^{(1)}(A)\in X_A$.

\begin{proof}
The spectrum of matrix  $A-{xx}^\top$ is equal to the spectrum of matrix $A$ up to
replacing  $\lambda_1(A)$ with zero. Since $X_A$ is not empty, from Lemma
\ref{lemma_lambda2_empty} if follows that all other eigenvalues are non-positive,
and so, matrix $A-{xx}^\top$ is negative semidefinite.
\end{proof}
\end{lemma}

\begin{lemma}\label{lemma_grow_x}
If ${x}\in X_A$, then $a{x}\in X_A$ for all $a>1$.

\begin{proof}
The proof is straightforward.
\end{proof}
\end{lemma}

\begin{lemma}\label{lemma_analytical_eq}
If $X_A\neq\emptyset$ and $X_A\neq\mathbb{R}^n$, then
\begin{equation}\label{eq_analytical_ineq}X_A =
\left\{{x}\in \mathbb{R}^n: \sum_{i=1}^n\frac{\left({x}^\top
u^{(i)}(A)\right)^2}{\lambda_i(A)}\ge 1\right\}.
\end{equation}

\begin{proof}
Let us denote with $I$ the identity $n\times n$ matrix. By definition of $X_A$, from
${x}\in X_A$ it follows that the characteristic equation $\det(A-{xx}^\top-\rho
I)=0$ has no positive roots. Since eigenvalues are continuous with respect to matrix
elements, identity $\lambda_1(A-{xx}^\top)=0$ holds on the boundary of $X_A$.
Therefore, if vector $x$ belongs to the boundary of $X_A$, then $\rho=0$ is a root
of the characteristic equation, i.e,
\begin{equation}\label{eq_boundary}
\det(A-{xx}^\top)=0.
\end{equation}

To solve equation (\ref{eq_boundary}), let us consider the spectral decomposition $U
\diag(\lambda) U^\top$ of matrix $A$, where $\lambda=(\lambda_i(A))_{i=1}^n$,
$U=(u^{(1)}(A), ..., u^{(n)}(A))$.

The characteristic equation and its roots are insensitive to orthogonal
transformations. Hence, $\det(A-{xx}^\top)=\det(\diag(\lambda)-{yy}^\top)$, where
${y}:=U^\top{x}$. Therefore, equation (\ref{eq_boundary}) can be written as
$$\det\left(%
\begin{array}{ccccc}
  \lambda_1(A)-y_1^2 & -y_1y_2 & \ldots & &-y_1y_n  \\
  -y_2y_1 & \lambda_2(A)-y_2^2 & -y_2y_3 & \ldots & -y_2y_n \\
  \vdots & -y_3y_2 & \ddots & &\vdots \\
  \vdots & \vdots &  & \ddots &\vdots \\
  -y_ny_1 & -y_ny_2 & & \ldots & \lambda_n(A)-y_n^2 \\
\end{array}%
\right)=0.$$ Let us transform the matrix to the triangular form with elementary row
operations not affecting the roots of the equation.

First we assume that $y_i\neq 0$, $i=1,...,n$. Let us divide $i$-th row by $y_i$,
$i=1,...,n$, and subtract the first row from all other rows obtaining the equation
$$\det\left(%
\begin{array}{ccccc}
  \frac{\lambda_1(A)}{y_1}-y_1 & -y_2 & \ldots & & -y_n  \\
  -\frac{\lambda_1(A)}{y_1} & \frac{\lambda_2(A)}{y_2} & 0 & \ldots & 0 \\
  \vdots & 0 & \ddots & & \vdots \\
  \vdots & \vdots &  & \ddots& \vdots \\
  -\frac{\lambda_1(A)}{y_1} & 0 & \ldots & &\frac{\lambda_n(A)}{y_n} \\
\end{array}%
\right)=0.$$

Let us multiply $i$-th row, $i=1,...,n$, by $\frac{y_i}{\lambda_i(A)}$ and add to
the first row all other rows, multiplying them by $y_i$. Finally we obtain the
desired lower triangular form:
\begin{equation}\label{eq_triang}\det\left(%
\begin{array}{ccccc}
  \frac{\lambda_1(A)}{y_1}-y_1 - \sum_{i=2}^n-\frac{y_i^2\lambda_1(A)}{y_1\lambda_i(A)}& 0 & &\ldots & 0  \\
  -\frac{y_2\lambda_1(A)}{y_1\lambda_2(A)} & 1 & 0 & \ldots & 0 \\
  \vdots & 0 & \ddots &  & \vdots \\
  \vdots & \vdots &  & \ddots & \vdots \\
  -\frac{y_n\lambda_1(A)}{y_1\lambda_n(A)} & 0 & \ldots & & 1 \\
\end{array}%
\right)=0.
\end{equation}

The determinant of a triangular matrix is equal to the product of its diagonal
element, so, equation (\ref{eq_triang}) can be written as
$$\frac{\lambda_1(A)}{y_1}-y_1 - \sum_{i=2}^n\frac{y_i^2\lambda_1(A)}{y_1\lambda_i(A)}=0.$$

Multiplying both sides of the equation by $\frac{y_1}{\lambda_1(A)}$, we finally
obtain
\begin{equation}\label{eq_hyperboloid}
\sum_{i=1}^n\frac{y_i^2}{\lambda_i(A)}= 1.
\end{equation}

If $y_i=0$ for some $i$, $i$-th row is already diagonal, no transformation needed,
so the case when some $y_i$ are equal to zero is considered is a similar manner.

Since $X_A$ is not empty and $X_A\neq\mathbb{R}^n$, it follows from Lemmas
\ref{lemma_lambda2_empty} and \ref{lemma_negative_semidefinite} that
$\lambda_1(A)>0$, $\lambda_2(A)\le0$. Therefore, equation (\ref{eq_hyperboloid})
defines the two-sheet hyperboloid in the $n$-dimensional space:
\begin{equation}\label{eq_hyperboloid_bisheet}
\frac{y_1^2}{\lambda_1(A)}-\sum_{i=2}^n\frac{y_i^2}{|\lambda_i(A)|}= 1.
\end{equation}
The boundary of the set $X_A$ belongs to this hyperboloid. Using Lemmas
\ref{lemma_perron_vector} and \ref{lemma_grow_x} one can easily check that both
sheets defined by the inequality
\begin{equation}
\frac{y_1^2}{\lambda_1(A)}-\sum_{i=2}^n\frac{y_i^2}{|\lambda_i(A)|}\ge 1
\end{equation}
have points from $X_A$ and, hence, belong to $X_A$. The space between these sheets
does not belong to $X_A$, since point ${x}={u}_2(A)$ is obviously does not belong to
$X_A$.

Taking into account that ${y}=U^\top {x}$, we obtain the desired inequality.
\end{proof}
\end{lemma}


\end{document}